\documentclass{amsart}\usepackage{amssymb,graphicx}

\theoremstyle{plain}
\newtheorem{Thm}{Theorem}
\newtheorem{Prop}[Thm]{Proposition}
\newtheorem{Lem}[Thm]{Lemma}

\newcommand{\z}{\textstyle}
\newcommand{\iE}{{\rm E}}\newcommand{\iG}{{\rm G}}\newcommand{\iM}{{\rm M}}\newcommand{\iS}{{\rm S}}

\newcommand{\C}{\mathbb{C}}\newcommand{\M}{\mathbb{M}}\newcommand{\N}{\mathbb{N}}
\newcommand{\Q}{\mathbb{Q}}\newcommand{\R}{\mathbb{R}}\newcommand{\Z}{\mathbb{Z}}
\newcommand{\cH}{\mathcal{H}}\newcommand{\cM}{\mathcal{M}}\newcommand{\fM}{\mathfrak{M}}
\newcommand{\cO}{\mathcal{O}}\newcommand{\cS}{\mathcal{S}}

\newcommand{\BM}{\begin{smallmatrix}}\newcommand{\EM}{\end{smallmatrix}}
\newcommand{\BC}{\begin{cases}}\newcommand{\EC}{\end{cases}}
\newcommand{\VS}[1]{\z\sum\limits_{#1}}\newcommand{\VP}[1]{\z\prod\limits_{#1}}
\newcommand{\VT}[1]{\z\bigoplus\limits_{#1}}
\newcommand{\6}{\;\;\;\;\;\;}
\newcommand{\ang}[1]{\langle #1\rangle}\newcommand{\TL}{\vartriangleleft}
\newcommand{\hto}{\hookrightarrow}\newcommand{\bto}{\xrightarrow{\sim}}

\begin{document}
\title{explicit structure of the graded ring of modular forms}
\author{Suda Tomohiko}
\email{t@mshk1201.com}
\maketitle

\section{Introduction}

For each $k\in\M=\{0,1,2,\cdots\}$ and a congruence subgroup $\Gamma\subset{\rm SL}_2(\Z)$, let $\cM(\Gamma)_k$ be the $\C$-vector space of all modular forms of weight $k$ with respect to $\Gamma$ and $\cS(\Gamma)_k$ be the subspace of all cusp forms. We study the $\M$-graded ring
\[\cM(\Gamma)=\VT{k\in\M}\cM(\Gamma)_k\]
and the homogeneous ideal
\[\cS(\Gamma)=\VT{k\in\M}\cS(\Gamma)_k.\]
Note $\cM(\Gamma)_0=\C$ and $\cS(\Gamma)_0=\{0\}$. We abbreviate
\[\cM(N)=\cM(\Gamma(N))\]
\[\cS(N)=\cS(\Gamma(N))\]
where
\[\Gamma(N)=\big\{(\BM a&b\\c&d\EM)\in{\rm SL}_2(\Z)\,\big|\,(\BM a&b\\c&d\EM)\equiv(\BM 1&0\\0&1\EM)\mod N\big\}.\]

When $(\BM 1&n\\0&1\EM)\in\Gamma$, we regard $\cM(\Gamma)\subset\C[[q^{\frac1n}]]$ via the Fourier expansion, where
\[q^{\frac1n}=e^{\frac{2\pi iz}n}\6(z\in\cH).\]
Note $\cS(\Gamma)\subset\C[[q^{\frac1n}]]q^{\frac1n}$.

SAGE, a Software for Algebra and Geometry Experimentation, is completely free and packages together a wide range of open sorce mathmatics software. So, SAGE commands are stated for some results.

\newpage

{\bf Notation.}

$\N=\{1,2,\cdots\}$

$\M=\N\cup\{0\}$

$\Z$ The integres

$\Q$ The rational numbers

$\R$ The real numbers

$\C$ The complex numbers

$\zeta_n=e^{\frac{2\pi i}n}$ \big($\zeta_2=-1$, $\zeta_3=\omega=\frac{1-\sqrt3i}2$, $\zeta_4=i$\big)

$\cH=\{z\in\C \,|\, \Im z>0\}$

$\cO(X)=\{f:X\to\C \,|\, f:{\rm holomorphic}\}$

[ ] The Gauss symbol

$a^+=\BC a &\text{if }0<a\\0 &\text{oth}\EC$

$\delta_{i,j}=\BC 1 &\text{if }i=j\\0 &\text{oth}\EC$\\

$R[x]$ The polynomial ring with a ring $R$ and a variable $x$

$R[[x]]$ The formal power series ring

$\natural$ natural map

$\ang{x_1,x_2,\cdots,x_n}$ the subgroup generated by $x_1,x_2,\cdots,x_n$

\newpage

\section{For the case $N=1$}

\subsection{Modular form of level 1}

First note $\Gamma(1)={\rm SL}_2(\Z)$. Put
\[{\rm GL}_2^+(\Q)=\big\{(\BM a&b\\c&d\EM) \,\big|\, a,b,c,d\in\Q,\,ad-bc>0\big\},\]
and define the action ${\rm GL}_2^+(\Q)\curvearrowright\cH$ by $(\BM a&b\\c&d\EM)z=\z\frac{az+b}{cz+d}$.

Moreover, for $k\in\M$ let $|_k:\cO(\cH)\curvearrowleft{\rm GL}_2^+(\Q)$ by
\[f|_k(\BM a&b\\c&d\EM):z\mapsto\frac{f((\BM a&b\\c&d\EM)z)}{(cz+d)^k}\]
then
\[\cM(1)_k=\big\{f\in\cO(\cH) \,\big|\, \forall\gamma\in\Gamma(1).f|_k\gamma=f,\,f:{\rm holomorphic\;at\;}\infty\big\}.\]
The word "$f:{\rm holomorphic\;at\;}\infty$" means that $\lim_{\Im z\to\infty}f(z)$ exists or even just that $f(z)$ is bounded as $\Im z\to\infty$.\\

Eisenstein series are representative modular forms. For $k\in2\N$, put
\[\iE_{1,k}(z)=\z\frac1{2\zeta(k)}\VS{(a,b)\in(\Z\times\Z)\setminus\{(0,0)\}}\frac1{(az+b)^k}=1-\frac{2k}{B_k}\VS{n\in\N}\VS{d|n}d^{k-1}q^n\]
where $\zeta$ is the Riemann zeta function and $B_k$ is the $k$-th Bernoulli number,
\begin{align*}\iE_{1,2}&=1-24\VS{n\in\N}\VS{d|n}dq^n\\
\iE_{1,4}&=1+240\VS{n\in\N}\VS{d|n}d^3q^n\\
\iE_{1,6}&=1-504\VS{n\in\N}\VS{d|n}d^5q^n\\
\iE_{1,8}&=1+480\VS{n\in\N}\VS{d|n}d^7q^n\end{align*}
for example. It is well-known $\iE_{1,k}\in\cM(1)_k$ for $k\geq4$.

As for $k=2$, we see for $(\BM a&b\\c&d\EM)\in\Gamma(1)$
\[\iE_{1,2}|_2(\BM a&b\\c&d\EM)(z)=\iE_{1,2}(z)-\frac{6ci}{\pi(cz+d)}.\]
$\iE_{1,2}-\frac3\pi{\rm Im}$ is $|_2$ invariant under $\Gamma(1)$, but it is not holomorphic (\cite[p18]{DS}).

We make $\iE_{1,4}+O(q^{20})$ in SAGE by

\verb~prec=20;eisenstein_series_qexp(4,prec,normalization='constant')~
\\

The dimension formula is also well-known : for $k\in\M$
\[\dim\cM(1)_{2k}=\big[\z\frac k6\big]+\BC0&\text{if }k\in6\M+1\\1&\text{otherwise}\EC\]

In particular $\dim\cM(1)_8=1$ and $\iE_{1,8}-\iE_{1,4}^2\in\cM(1)_8\cap\C[[q]]q=\{0\}$. That is $\iE_{1,8}=\iE_{1,4}^2$ and a relation between divisor sums is obtained : for $n\in\N$
\[\VS{d|n}d^7=\VS{d|n}d^3+120\VS{i=1}^{n-1}\Big(\VS{d|i}d^3\VS{d|(n-i)}d^3\Big).\]
Similarly $\iE_{1,10}=\iE_{1,4}\iE_{1,6}$ and $\iE_{1,14}=\iE_{1,4}^2\iE_{1,6}$.

\newpage

\subsection{Dedekind eta function}

Put
\[\eta=q^{\frac1{24}}\VP{n\in\N}(1-q^n)\]
then $\eta\in\cO(\cH)$ and $\eta$ is non-zero on $\cH$. Euler's pentagonal numbers theorem says
\[\VP{n\in\N}(1-q^n)=\VS{n\in\Z}(-1)^nq^{\frac{3n^2-n}2}\]
thus
\[\eta=\VS{n\in\Z}(-1)^nq^{\frac{(6n-1)^2}{24}}=\VS{n\in\N}(\frac 3n)q^{\frac{n^2}{24}}\]
where $(\frac nm)$ is the Kronecker symbol.\\

Now, for $(\BM a&b\\c&d\EM)\in{\rm SL}_2(\Z)$ and $\kappa\in\frac12\M$ let
\[\z f|_\kappa(\BM a&b\\c&d\EM):z\mapsto ((\frac cd)\sqrt{cz+d})^{-2\kappa}f((\BM a&b\\c&d\EM)z)\]
where $\sqrt x$ is the "principal" determination of the square root of $x$ i.e. the one with argument in $(-\frac\pi2,\frac\pi2]$. It is not an action, however satisfies $f^2|_1\gamma=(f|_{\frac12}\gamma)^2$.

In fact, we can't make weight $\frac12$ action properly. Indeed $f|_{\frac12}(\BM -1&0\\0&-1 \EM)$ should be $\pm if$ since $f^2|_1(\BM -1&0\\0&-1 \EM)=-f^2$, contradict to $f|_{\frac12}(\BM -1&0\\0&-1 \EM)|_{\frac12}(\BM -1&0\\0&-1 \EM)=f$.\\

Note the logarithmic derivative
\[\frac d{dz}\log\eta(z)=\frac{\pi i}{12}+2\pi i\VS{n\in\N}\dfrac{nq^n}{1-q^n}=\dfrac{\pi i}{12}\iE_{1,2}(z).\]

The eta function satisfies the functional equations
\[\z\eta|_{\frac12}(\BM 1&1\\0&1 \EM)=e^{\frac{\pi i}{12}}\eta,\6\eta|_{\frac12}(\BM 0&-1\\1&0 \EM)=e^{-\frac{\pi i}{4}}\eta.\]
So $\eta$ is essentially a $\frac12$-weight form of ${\rm SL}_2(\Z)=\ang{(\BM 1&1\\0&1 \EM),(\BM 0&-1\\1&0 \EM)}$ and $\eta^{24}\in\cM(1)_{12}$,
\[\z\frac1{12^3}(\iE_{1,4}^3-\iE_{1,6}^2)=\eta^{24}.\]

We make $\eta^{24}+O(q^{20})$ in SAGE by

\verb~R.<q>=PowerSeriesRing(QQ);prec=20~

\verb~from sage.modular.etaproducts import qexp_eta~

\verb~E=qexp_eta(ZZ[['q']],prec);E^24*q~\\[5pt]

Moreover $\eta^{24}\in\cS(1)_{12}$ since $\Gamma(1)$ has only one cusp $i\infty$ and
\[\cM(1)\cap\C[[q]]q=\cS(1).\]

Rademacher showed for $(\BM a&b\\c&d \EM)\in\Gamma(1)$
\[\eta|_{\frac12}(\BM a&b\\c&d \EM)=\z e^{\frac{\pi i}{12}(ac+bd-cd-acd^2+3d-3)}\eta\6\text{if }c\in2\N,\]
\[\eta|_{\frac12}(\BM a&b\\c&d \EM)=\z(\frac dc)(\frac cd)e^{\frac{\pi i}{12}(ac+bd+cd-bc^2d-3c)}\eta\6\text{if }c\in2\M+1.\]

\newpage

\subsection{Graded ring}

Let $S$ be a commutative semigroup and $R=\VT{s\in S}R_s$ is a $S$-graded ring, i.e.
\[R=\VT{s\in S}R_s,\6R_s\cdot R_t\subset R_{s+t}\text{ for }s,t\in S.\]
If $\dim R_s<\infty$ for all $s\in S$, then define the Hilbert function
\[{\rm Dim}\,R=\VS{s\in S}(\dim R_s)t^s.\]

For example, $\C$ is $S$-graded by $\C_0=\C$ and $\C_s=\{0\}$ for $s\in S\setminus\{0\}$ where 0 is the unit element of $S$. We regard $1=t^0$, then ${\rm Dim}(\C)=1$.

For a subgraoup $T\subset S$, we define $R|_T=\VT{t\in T}R_t$ and $\VS{s\in S}a_st^s\big|_T=\VS{s\in T}a_st^s$, then ${\rm Dim}(R|_T)=({\rm Dim}\,R)|_T$.\\

For $R=\VT{s\in S}R_s$ and $s_1,\cdots,s_r\in\N$, define
\[R'=R[X_1,\cdots,X_r]^{[s_1,\cdots,s_r]}\]
to be the ring $R[X_1,\cdots,X_r]$ graded by $R_s\subset R'_s$ and $X_i\in R'_{s_i}$.
\begin{Lem}\label{HF}${\rm Dim}\big(\C[X,Y]^{[m,n]}\big)=\dfrac1{(1-t^m)(1-t^n)}$.\end{Lem}
\begin{proof}$(R[X]^{[n]})_k=R_k\oplus R_{k-n}X\oplus R_{k-2n}X^2\oplus\cdots$ and ${\rm Dim}\big(R[X]^{[n]}\big)=\dfrac{{\rm Dim}(R)}{1-t^n}$.\end{proof}

We calculate
\[\frac1{(1-t)(1-t^n)}=\VS{k\in\M}\big([\frac kn]+1\big)t^k\]
and
\begin{align*}\frac1{(1-t^2)(1-t^3)}&=\frac1{t^2}\Big(\frac1{(1-t)(1-t^2)}-\frac1{(1-t)(1-t^3)}\Big)\\
&=\VS{k\in\M}\big([\frac k2]-[\frac k3]\big)t^{k-2}\end{align*}

\begin{Thm}$\natural:\C[\iE_{1,4},\iE_{1,6}]^{[4,6]}\bto\cM(1)$ and $\cM(1)\eta^{24}=\cS(1)$.

\end{Thm}
\begin{proof}Put $R=\C[\iE_{1,4},\iE_{1,4}^3-\iE_{1,6}^2]$ then $\natural:R\hto\cM(1)|_{4\M}$.

Since $\C[\iE_{1,4},\iE_{1,6}]=R\oplus R\iE_{1,6}$, we get $\hto$ of the former assertion.

It follows from ${\rm Dim}\,\C[\iE_{1,4},\iE_{1,6}]^{[4,6]}={\rm Dim}\,\cM(1)$.

The latter one follows from $\cM(1)\eta^{24}\subset\cS(1)$ and
\[\dim\cS(1)_k=(\dim\cM(1)_k-1)^+=\dim\cM(1)_{k-12}\]
for $k\in\M$.
\end{proof}

\newpage

\subsection{Ring with square}

For a ring $R$ and $\phi\in R$, $n\in\M+2$ let
\[R\big[\sqrt[n]\phi\big]=R[X]\big/(X^n-\phi).\]
The map $\sqrt[3]1\mapsto\frac{-1+\sqrt{-3}}2:\Q\big[\sqrt[3]1\big]\to\Q\big[\sqrt{-3}\big]$ is well-defined (and not injective).

We abbreviate $R\big[\sqrt{X},\sqrt{Y}\big]=\big(R\big[\sqrt{X}\big]\big)\big[\sqrt{Y}\big]$ and so on. Remark that the natural map
\[\natural:\C\big[X,\sqrt{XY},\sqrt{YZ},\sqrt{ZX}\big]\to\C\big[\sqrt{X},\sqrt{Y},\sqrt{Z}\big]\]
is not injective since $\natural\big(X\cdot\sqrt{YZ}\big)=\natural\big(\sqrt{XY}\cdot\sqrt{ZX}\big)$.\\

If $R$ is a $\M$-graded ring and $\phi\in R_{nk}\setminus R{\text \textasciicircum}n$ where $A{\text \textasciicircum}n=\{a^n \,|\, a\in A\}$, then we can make the $\M$-graded ring
\[R\big[\sqrt[n]\phi\big]^{[k]}.\]
\begin{Lem}${\rm Dim}\big(R\big[\sqrt[n]\phi\big]^{[k]}\big)={\rm Dim}(R)\dfrac{1-t^{nk}}{1-t^k}$.\end{Lem}
\begin{proof}$R[\sqrt\phi]=R\oplus R\sqrt[n]\phi\oplus R\sqrt[n]\phi^2\cdots\oplus R\sqrt[n]\phi^{n-1}$.\end{proof}

If $\phi\in\big(\C[X,Y]^{[1,1]}\big)_2\setminus\C[X,Y]{\text \textasciicircum}2$ then
\[{\rm Dim}\big(\C[X,Y,\z\sqrt\phi]^{[1,1,1]}\big)=\dfrac{1+t}{(1-t)^2}=\VS{k\in\M}(2k+1)t^k.\]\\

Deligne's theorem says
\[\natural:\Z\big[\iE_{1,4},\eta^{24},\iE_{1,6}\big]^{[4,12,6]}\bto\cM(1)\cap\Z[[q]]\]
where $\iE_{1,6}=\sqrt{\iE_{1,4}^3-12^3\eta^{24}}$.\\

We have another meaning of $\sqrt{\;}$. Generally, for $\phi\in q^{rn}+\C[[q]]q^{rn+1}$, we define the $n$-th square $\sqrt[n]\phi\in q^r+\C[[q]]q^{r+1}$ by $(\sqrt[n]\phi)^n=\phi$. For example $\sqrt{q}=q^{\frac12}$,
\[\sqrt{1+q}=\z 1+\frac12q-\frac18q^2+\frac1{16}q^3-\frac5{128}q^4+\cdots.\]
Remark for a modular form $f\in\C[[q]]$, in general $\sqrt f(z)\neq\sqrt{f(z)}$.

SAGE's command \verb~sqrt~ return $\sqrt{\;}$ by Newton's method. We can make $\sqrt[3]{\;}$ by

\verb~R.<q>=PowerSeriesRing(QQ);prec=20~

\verb~def sqrt3(x):~

\verb~    j=x.exponents()[0]~

\verb~    if j%3>0 or x.list()[j]<>1:~

\verb~        return "Error"~

\verb~    j=j/3;s=1~

\verb~    for i in range(1,prec-j):~

\verb~        ary=(x/q^(3*j)-s^3).list()~

\verb~        if i<=len(ary):~

\verb~            s=s+ary[i]/3*q^i~

\verb~    return s*q^j+O(q^prec)~

\newpage

\section{For the case $N=2$}

\subsection{Congruence group}

For $N\in\N$, define
\[\Gamma_0(N)=\big\{(\BM a&b\\c&d\EM)\in{\rm SL}_2(\Z) \,\big|\, c\in N\Z\big\},\]
\[\Gamma^0(N)=\big\{(\BM a&b\\c&d\EM)\in{\rm SL}_2(\Z) \,\big|\, b\in N\Z\big\},\]
\[\Gamma_1(N)=\big\{(\BM a&b\\c&d\EM)\in\Gamma_0(N) \,\big|\, d\in N\Z+1\big\},\]
\[\Gamma^1(N)=\big\{(\BM a&b\\c&d\EM)\in\Gamma^0(N) \,\big|\, d\in N\Z+1\big\},\]
\[\Gamma(N)=\Gamma_1(N)\cap\Gamma^1(N).\]

Let the natural map $\natural_N:\Z\to\Z/N$, then
\[(\BM a&b\\c&d\EM)\mapsto(\BM \natural_N(a)&\natural_N(b)\\\natural_N(c)&\natural_N(d)\EM):{\rm SL}_2(\Z)/\Gamma(N)\bto{\rm SL}_2(\Z/N),\]
and $[{\rm SL}_2(\Z):\Gamma(N)]=N^3\VP{p|N}(1-\frac1{p^2})$. SAGE command for $N=12$ :

\verb~N=12~

\verb~Gamma(N).index()-N^3*prod([1-1/f[0]^2 for f in list(factor(N))])~

$\longrightarrow$ \verb~0~\\

A subgroup $H\subset{\rm SL}_2(\Z)$ is called a congruence subgroup if there exists $N$ such that $\Gamma(N)\subset H$. The level of $H$ is then the smallest such $N$.\\

Generally, for groups $G$ and $H$, denote the group of all homomorphisms $G\to H$ by ${\rm Hom}(G,H)$. We denote the trivial character in ${\rm Hom}(G,\C^\times)$ by ${\tt 1}_G$.

Let
\[d_N:(\BM a&b\\c&d\EM)\mapsto\natural_N(d),\]
then $d_N\in{\rm Hom}(\Gamma_0(N),\Z/N^\times)$. When $d_N\in{\rm Hom}(\Gamma,\Z/N^\times)$ for a congruence subgroup $\Gamma$, for $\chi\in{\rm Hom}(\Z/N^\times,\C^\times)$ we regard $\chi=\chi\circ d_N\in{\rm Hom}(\Gamma,\C^\times)$.

For a subgroup $G\subset\Z/N^\times$, we define
\[\Gamma_G(N)=\{\gamma\in\Gamma_0(N) \,|\, d_N(\gamma)\in G\}.\]
In particular, $\Gamma_{\Z/N^\times}(N)=\Gamma_0(N)$ and $\Gamma_{\{1\}}(N)=\Gamma_1(N)$. We have
\[d_N:\Gamma_G(N)/\Gamma_1(N)\bto G.\]\\

Generally, for a group $G$ we put an action ${\TL}:G\curvearrowleft G$ by $g\TL h=h^{-1}gh$. Since
\[(\BM a&b\\c&d\EM)\TL(\BM h&0\\0&1\EM)=(\BM a&b/h\\ch&d\EM)\]
we see
\[\Gamma_{\ang{N+1}}(N^2)=\Gamma(N)\TL(\BM N&0\\0&1\EM),\]
for example $\Gamma_0(4)=\Gamma(2)\TL(\BM 2&0\\0&1\EM)$ and $\Gamma_{\ang4}(9)=\Gamma(3)\TL(\BM 3&0\\0&1\EM)$.

\newpage

\subsection{Modular form}

For a congruence subgroup $\Gamma\subset{\rm SL}_2(\Z)$, we define
\[\cM(\Gamma)_k=\big\{f\in\cO(\cH) \,\big|\, \forall\gamma\in\Gamma.f|_k\gamma=f,\,\forall\alpha\in{\rm SL}_2(\Z).f|_k\alpha:{\rm holomorphic\;at\;}\infty\big\}.\]

If $\Gamma\subset\Gamma'$ then $\cM(\Gamma)\supset\cM(\Gamma')$.

If $(\BM-1&0\\0&-1\EM)\in\Gamma$ and $k\in2\M+1$ then $\cM(\Gamma)_k=\{0\}$.

For $\alpha\in{\rm GL}_2^+(\Q)$ we may construct the isomorphism of $\M$-graded rings
\[|\alpha=\VT{k\in\M}|_k\alpha:\cM(\Gamma)\bto\cM(\Gamma\TL\alpha).\]\\

Define $f^{\ang h}=f|(\BM h&0\\0&1\EM)$, then $q^{\ang h}=q^h$.

\begin{Lem}\label{V1}For $N,h\in\N$ we have
\[\cM(\Gamma_0(N))\subset\cM(\Gamma_0(hN)),\]
\[\cM(\Gamma_0(N))^{\ang h}=\cM(\Gamma_0(hN))\cap\C[[q^h]].\]
\end{Lem}
\begin{proof}$\Gamma_0(N)\supset\Gamma_0(hN)$. $\Gamma_0(N)\TL(\BM h&0\\0&1\EM)$ is generated by $\Gamma_0(hN)$ and $(\BM 1&\frac1h\\0&1\EM)$.\end{proof}

We see $q^{\frac1N}|(\BM 1&1\\0&1\EM)=\zeta_Nq^{\frac1N}$. Note $\Gamma(N)\TL(\BM 1&1\\0&1\EM)=\Gamma(N)$ since $\Gamma(N)$ is a normal subgroup of ${\rm SL}_2(\Z)$, hence
\[\cM(N)|(\BM 1&1\\0&1\EM)=\cM(N).\]

For $N\in\M+2$, let the Eisenstein series
\[\iE_{N,2}=\z\frac1{N-1}\big(N\iE_{1,2}^{\ang N}-\iE_{1,2}\big)=1+\frac{24}{N-1}\VS{n\in\N}\VS{d|n,N\nmid d}dq^n\]
then $\iE_{N,2}\in\cM(\Gamma_0(N))_2$.

\begin{Lem}\label{cuspform}$\eta^{12}\in\cS(2)_6$, $\eta^8\in\cS(3)_4$, $\eta^6\in\cS(4)_3$, $\eta^4\in\cS(6)_2$ and $\eta^2\in\cS(12)_1$.\end{Lem}
\begin{proof}For $N=2,3,4,6,12$ and $(\BM a&b\\c&d \EM)\in\Gamma(1)$ we see
\[\eta^{\frac{24}N}|_{\frac{12}N}(\BM a&b\\c&d \EM)=e^{\frac{2\pi i}N(bd+3d-3)}\eta^{\frac{24}N}\6\text{if }c\in2\N\cap N\N,\]
As for $N=3$,
\[\eta^8|_4(\BM a&b\\c&d \EM)=e^{\frac{2\pi i}3bd}\eta^8\6\text{if }\z c\in(2\M+1)\cap 3\N.\]\end{proof}

For a congruence group $\Gamma$ and $\chi\in{\rm Hom}(\Gamma,\C^\times)$, we define
\[\iM(\Gamma)_{\chi,k}=\big\{f\in\cM(\ker\chi)_k \,\big|\, f|_k\gamma=\chi(\gamma)f\text{ for all }\gamma\in\Gamma\big\}.\]
We see $\iM(\Gamma)_{{\tt 1}_\Gamma,k}=\cM(\Gamma)_k$ and
\[(\BM -1&0\\0&-1\EM)\in\Gamma,\;\chi(\BM -1&0\\0&-1\EM)\neq(-1)^k\Longrightarrow\iM(\Gamma)_{\chi,k}=\{0\},\]
\[\iM(\Gamma)_{\chi,k}\iM(\Gamma)_{\chi',k'}\subset\iM(\Gamma)_{\chi\chi',k+k'}.\]

We abbreviate $\iM(N)=\iM(\Gamma_0(N)\cap\Gamma^0(N))$.

\newpage

\subsection{Dimension formula}

See \cite[\S3]{DS}, \cite[\S6.3]{St}.

For a congruence group $\Gamma$ of ${\rm SL}_2(\Z)$, let
\[X(\Gamma)=\Gamma\setminus(\cH\cup\Q\cup\{\infty\}),\]
which is a compact Riemann surface. Put
\[\varepsilon_2(\Gamma)=\text{the number of elliptic points of period 2 in }X(\Gamma)\]
\[\varepsilon_3(\Gamma)=\text{the number of elliptic points of period 3 in }X(\Gamma)\]
\[\varepsilon_\infty(\Gamma)=\text{the number of cusps of }X(\Gamma)\]
\[d(\Gamma)={\rm deg}\big(X(\Gamma)\to X({\rm SL_2}(\Z))\big)=\BC [{\rm SL}_2(\Z):\Gamma]&{\rm if}\;(\BM-1&0\\0&-1\EM)\in\Gamma \\ \frac12[{\rm SL}_2(\Z):\Gamma] & {\rm oth}\EC\]
\[g(\Gamma)=\text{the genus of }X(\Gamma)=1+\frac{d(\Gamma)}{12}-\frac{\varepsilon_2(\Gamma)}4-\frac{\varepsilon_3(\Gamma)}3-\frac{\varepsilon_\infty(\Gamma)}2.\]
Then, for $k\in2\N$ we have
\[\dim\cM(\Gamma)_k=(k-1)(g(\Gamma)-1)+\Big[\frac k4\Big]\varepsilon_2(\Gamma)+\Big[\frac k3\Big]\varepsilon_3(\Gamma)+\frac k2\varepsilon_\infty(\Gamma)\]
\[\dim\cS(\Gamma)_k=\BC \dim\cM(\Gamma)_k-\varepsilon_\infty(\Gamma) & k\geqq4 \\ g(\Gamma) & k=2 \EC\]\\

In addtion, put
\[\varepsilon_\infty^{\rm reg}(\Gamma)=\text{the number of regular cusps of }X(\Gamma)\]
\[\varepsilon_\infty^{\rm irr}(\Gamma)=\text{the number of irregular cusps of }X(\Gamma)\]
Note $\varepsilon_\infty(\Gamma)=\varepsilon_\infty^{\rm reg}(\Gamma)+\varepsilon_\infty^{\rm irr}(\Gamma)$.

When $(\BM -1&0\\0&-1\EM)\notin\Gamma$, for $k\in2\N+1$ we have
\[\dim\cM(\Gamma)_k=(k-1)(g(\Gamma)-1)+\Big[\frac k3\Big]\varepsilon_3(\Gamma)+\frac k2\varepsilon_\infty^{\rm reg}(\Gamma)+\frac{k-1}2\varepsilon_\infty^{\rm irr}(\Gamma)\]
\[\dim\cS(\Gamma)_k=\dim\cM(\Gamma)_k-\varepsilon_\infty^{\rm reg}(\Gamma).\]

As for $k=1$ we have
\[\dim\cM(\Gamma)_1\BC =\frac 12\varepsilon_\infty^{\rm reg}(\Gamma) & \text{if }\varepsilon_\infty^{\rm reg}(\Gamma)>2g(\Gamma)-2 \\[0.1cm] \geq\frac 12\varepsilon_\infty^{\rm reg}(\Gamma) & \text{if }\varepsilon_\infty^{\rm reg}(\Gamma)\leq2g(\Gamma)-2 \EC\]
\[\z\dim\cS(\Gamma)_1=\dim\cM(\Gamma)_1-\frac 12\varepsilon_\infty^{\rm reg}(\Gamma)\]\\

Now, we prepare datas for $\Gamma(N)$ : for $N\in\M+2$
\[\varepsilon_\infty^{\rm reg}(\Gamma(N))=d(\Gamma(N))/N\]
\[\varepsilon_2(\Gamma(N))=\varepsilon_3(\Gamma(N))=\varepsilon_\infty^{\rm irr}(\Gamma(N))=0\]\\

SAGE command for $\dim\cM(\Gamma_{\ang4}(9))_k$ $(1\le k\le10)$ :

\verb~[GammaH(9,[4]).dimension_modular_forms(k) for k in range (1,11)]~

$\longrightarrow$ \verb~[2, 3, 4, 5, 6, 7, 8, 9, 10, 11]~

\newpage

\subsection{$/\;/$ and $"\;"$ operator}

For a prime number $p$, put
\[f^{/p/}=\frac{f^{\ang{\frac1p}p}}f,\6 f^{"p"}=\frac{f^{\ang pp}}f\]

We see $\eta^{\ang h}\big((\BM 0&-1\\1&0\EM)z\big)=\sqrt{\frac zh}e^{-\frac{\pi i}4}\eta^{\ang{\frac1h}}(z)$ and
\[\eta^{"p"}|_{\frac{p-1}2}(\BM 0&-1\\1&0\EM)=\frac{\eta^{\ang p}\big((\BM 0&-1\\1&0\EM)z\big)^p}{\sqrt z\eta\big((\BM 0&-1\\1&0\EM)z\big)}=\frac{e^{-\frac{(p-1)}4\pi i}}{\sqrt p^p}\eta^{/p/}.\]

Note
\[(\BM a&b\\c&d\EM)\TL(\BM 0&-1\\1&0\EM)=(\BM d&-c\\-b&a\EM).\]\\[0.5cm]

We abbreviate $\flat={/2/}$ and $\sharp="2"$.

\begin{Lem}$\eta^{\flat4}\in\iM(\Gamma^0(2))_{c_2,2}$ and $\eta^{\sharp4}\in\iM(\Gamma_0(2))_{b_2,2}$ where
\[b_2:(\BM a&b\\c&d\EM)\mapsto(-1)^b,\6c_2:(\BM a&b\\c&d\EM)\mapsto(-1)^c.\]
\end{Lem}
\begin{proof}For $(\BM a&b\\c&d \EM)\in\Gamma^0(2)$, we see
\[\frac{e^{\frac{2\pi i}3(2ac+\frac b2d-2cd-2acd^2+3d-3)}}{e^{\frac{\pi i}3(ac+bd-cd-acd^2+3d-3)}}=e^{\pi i(ac-cd-acd^2+d-1)}=(-1)^c,\]
\[\frac{e^{\frac{2\pi i}3(2ac+\frac b2d-2cd-2acd^2+3d-3)}}{e^{\frac{\pi i}3(ac+bd+cd-bc^2d-3c)}}=e^{\pi i(ac+\frac{-5-4ad+bc}3cd+d-1+c)}=(-1)^c.\]

Note $\eta^{\sharp}|_{\frac12}(\BM 0&-1\\1&0\EM)=\frac1{2e^{\frac{\pi i}4}}\eta^{\flat}$. The latter assertion  follows from the first one and $\eta^{\sharp4}|_2(\BM 0&-1\\1&0\EM)=-\frac1{16}\eta^{\flat4}$.
\end{proof}

\begin{Thm}$\natural:\C\big[\eta^{\flat4},\eta^{\sharp4}\big]^{[2,2]}\bto\cM(2)$ and $\cM(2)\eta^{12}=\cS(2)$.\end{Thm}
\begin{proof}The dimension formula says for $k\in\M$
\[\dim\cM(2)_k=k+1,\6\dim\cS(2)_k=(\dim\cM(2)_k-3)^+=\dim\cM(2)_{k-6}.\]
\end{proof}

\vspace{0.5cm}

For a prime number $p$, we see
\[f^{/p/"p"}=\frac{f^{/p/\ang pp}}{f^{/p/}}=\frac{f^{p^2+1}}{f^{\ang{\frac1p}p}f^{\ang{p}p}}=\frac{f^{"p"\ang{\frac1p}p}}{f^{"p"}}=f^{"p"/p/}.\]
Also, for another prime number $s$
\[f^{/p//s/}=\frac{f^{/p/\ang{\frac1s}s}}{f^{/p/}}=\frac{f^{\ang{\frac1{ps}}ps}f}{f^{\ang{\frac1s}s}f^{\ang{\frac1p}p}}=f^{/s//p/}\]
similarly $f^{/p/"s"}=f^{"s"/p/}$ and $f^{"p""s"}=f^{"s""p"}$.\\

We abbreviate $\natural=\flat\sharp$ i.e. $f^\natural=\dfrac{f^5}{f^{\ang{\frac12}2}f^{\ang22}}$. It is important that $f^{\flat}f^{\natural}=f^{\flat\ang22}$ and $f^{\sharp}f^{\natural}=f^{\sharp\ang{\frac12}2}$. We also abbreviate $f^{\Diamond}=f^{\flat}f^{\sharp}$ then $f^{\natural}f^{\Diamond}=f^3$.\\

\begin{Lem}\label{rel21}$\eta^{\natural4}=\eta^{\flat4}+(2\eta^{\sharp})^4$.\end{Lem}
\begin{proof}From $\big\{f\in\cM(2)_2 \,\big|\, f|_2(\BM 0&-1\\1&0\EM)=-f\big\}=\C\eta^{\natural4}$ and
\[\eta^{\natural}|_{\frac12}(\BM 0&-1\\1&0\EM)=\frac{\eta\big((\BM 0&-1\\1&0\EM)z\big)^5}{\sqrt z\eta^{\ang{\frac12}}\big((\BM 0&-1\\1&0\EM)z\big)^2\eta^{\ang2}\big((\BM 0&-1\\1&0\EM)z\big)^2}=e^{-\frac{\pi i}4}\eta^{\natural}.\]
\end{proof}

Note
\[\eta^{\flat}|(\BM 1&1\\0&1\EM)=\VP{n\in\frac12\N\setminus\N}(1+q^n)^2\VP{n\in\N}(1-q^n)=\frac{\VP{n\in\N\setminus2\N}(1-q^n)^3\VP{n\in2\N}(1-q^n)}{\VP{n\in\frac12\N\setminus\N}(1-q^n)^2}=\eta^{\natural}.\]\\

We see $\big\{f\in\cM(2)_2 \,\big|\, f|_2(\BM 1&1\\0&1\EM)=f\big\}=\C\iE_{2,2}$ and $\iE_{2,2}=\frac12(\eta^{\natural4}+\eta^{\flat4})$.

Since $\big\{f\in\cM(2)_k \,\big|\, f|_k(\BM 1&1\\0&1\EM)=f|_k(\BM 0&-1\\1&0\EM)=f\big\}=\cM(1)_k$, we see
\[\z\iE_{1,4}=\frac12(\eta^{\natural8}+\eta^{\flat8}+(2\eta^{\sharp})^8),\]
\[\z\iE_{1,8}=\frac12(\eta^{\natural16}+\eta^{\flat16}+(2\eta^{\sharp})^{16}).\]\\

Put $\iE_{2,2}^\gets=\iE_{2,2}^{\ang{\frac12}}\big|(\BM 1&1\\0&1\EM)$ then $\iE_{2,2}^{\gets}|_2(\BM 0&-1\\1&0\EM)=\iE_{2,2}^{\gets}$ since
\[(\BM \frac1n&0\\0&1\EM)(\BM 1&1\\0&1\EM)(\BM 0&-1\\1&0\EM)=(\BM 1&0\\n&1\EM)(\BM \frac1n&0\\0&1\EM)(\BM 1&-1\\0&1\EM).\]

We see $\big\{f\in\cM(2)_2 \,\big|\, f|_2(\BM 0&-1\\1&0\EM)=f\big\}=\C\iE_{2,2}^{\gets}$ and $\iE_{2,2}^{\gets}=\eta^{\flat4}-(2\eta^{\sharp})^4$. Acting $(\BM 1&1\\0&1\EM)$, we get $\iE_{2,2}^{\ang{\frac12}}=\eta^{\natural4}+(2\eta^{\sharp})^4$.

We also see $2\iE_{2,2}|(\BM 0&-1\\1&0\EM)=-\eta^{\natural4}-(2\eta^{\sharp})^4=-\iE_{2,2}^{\ang{\frac12}}$, thus
\[\z\iE_{1,4}=\frac16\big((2\iE_{2,2})^2+\iE_{2,2}^{\ang{\frac12}2}+\iE_{2,2}^{\gets2}\big),\]
\[\iE_{1,6}=\iE_{2,2}\iE_{2,2}^{\ang{\frac12}}\iE_{2,2}^{\gets},\]
\[\z\iE_{1,8}=\frac1{18}\big((2\iE_{2,2})^4+\iE_{2,2}^{\ang{\frac12}2}+\iE_{2,2}^{\gets2}\big).\]
SAGE command for $\iE_{1,6}$ identity :

\verb~R.<q>=PowerSeriesRing(QQ);prec=20~

\verb~E12=eisenstein_series_qexp(2,prec,normalization='constant')~

\verb~E16=eisenstein_series_qexp(6,prec,normalization='constant')~

\verb~E22=2*E12(q^2)-E12~

\verb~E16(q^2)-E22(q^2)*E22*E22(-q)~

$\longrightarrow$ \verb~O(q^20)~

\newpage

\section{For the case $N=4,8$}

\subsection{Nebentype modular form}

We abbreviate ${\tt 1}_N={\tt 1}_{\Z/N^\times}$. For example
\[{\rm Hom}(\Z/2^\times,\C^\times)=\{{\tt 1}_2\}.\]

We naturally regard ${\rm Hom}(\Z/N^\times,\C^\times)\subset{\rm Hom}(\Z/hN^\times,\C^\times)$
via tha natural map $\Z/hN^\times\to\Z/N^\times$. For exmaple ${\tt 1}_4={\tt 1}_2$.

Define $\rho_4=(\frac{-1}\bullet)$ and $\rho_8=(\frac\bullet 2)=(\frac2\bullet)$. For $n\ge4$, define $\chi_{2^n}\in{\rm Hom}(\Z/2^{n\times},\C^\times)$ by $\chi_{2^n}(5)=\zeta_{2^{n-2}}$ and $\chi_{2^n}(-1)=1$. For exmaple, $\chi_{16}(5)=i$ and $\chi_{16}^2=\rho_8$. Then
\begin{align*}
{\rm Hom}(\Z/4^\times,\C^\times)&=\ang{\rho_4}\\
{\rm Hom}(\Z/8^\times,\C^\times)&=\ang{\rho_4,\rho_8}\\
{\rm Hom}(\Z/16^\times,\C^\times)&=\ang{\rho_4,\chi_{16}}\\
{\rm Hom}(\Z/32^\times,\C^\times)&=\ang{\rho_4,\chi_{32}}
\end{align*}

For $\chi\in{\rm Hom}(\Z/N^\times,\C^\times)$, define $\overline\chi\in{\rm Hom}(\Z/N^\times,\C^\times)$ by $(\overline\chi)(x)=\overline{\chi(x)}$. For example, $\overline{\chi_{16}}=\chi_{16}^3$.\\[0.5cm]

Let $\Gamma$ be a congruence group and it satisfies $d_N\in{\rm Hom}(\Gamma,\Z/N^\times)$. For a subgroup $X\subset{\rm Hom}(\Z/N^\times,\C^\times)$ put
\[\iM(\Gamma)_{X,k}=\VT{\chi\in X}\iM(\Gamma)_{\chi,k}\]
and make the $(X\times\M)$-graded ring
\[\iM(\Gamma)_X=\VT{k\in\M}\iM(\Gamma)_{X,k}.\]

Then, we have decomposition (cf. \cite[\S4.3]{DS} or \cite[Proposition 9.2]{St})
\[\cM(\Gamma_1(N))=\iM(\Gamma_0(N))_{{\rm Hom}(\Z/N^\times,\C^\times)}.\]
Indeed, define $\ang d:\cM(\Gamma_1(N))_k\to\cM(\Gamma_1(N))_k$ by
\[\ang df=f|_k(\BM a& b\\c&\delta\EM)\]
for any $(\BM a& b\\c&\delta\EM)\in\Gamma_0(N)$ with $\delta\in N\Z+d$ and $\pi_\chi:\cM(\Gamma_1(N))_k\to\cM(\Gamma_1(N))_k$ by
\[\pi_\chi=\z\frac1{\phi(N)}\VS{d\in\Z/N^\times}\overline{\chi}(d)\ang d\]
then the decomposition is given by $\VS{\chi\in{\rm Hom}(\Z/N^\times,\C^\times)}\pi_\chi$.

\begin{Prop}\label{dec1}For a subgroup $G\subset\Z/N^\times$, we have
\[\cM(\Gamma_G(N))=\iM(\Gamma_0(N))_H,\]
where $H=\{f\in{\rm Hom}(\Z/N^\times,\C^\times) \,|\, f(G)=\{1\}\}$.
\end{Prop}
\begin{proof}If $\chi\in H$ then $\cM(\Gamma_G(N))_k\supset\iM(\Gamma_0(N))_{\chi,k}$.

If $\chi\in{\rm Hom}(\Z/N^\times,\C^\times)$ and $f(G)\neq\{1\}$, then we can choose $\alpha\in\Gamma_G(N)$ such that $\chi(\alpha)\neq1$, therefore $\iM(\Gamma_0(N))_{\chi,k}\cap\cM(\Gamma_G(N))_k=\{0\}$.
\end{proof}

It is important that $\cM(N)=\iM(N)_{{\rm Hom}(\Z/N^\times,\C^\times)}$.\\

Note
\[f\mapsto f^c:\iM(\Gamma_0(N))_{\chi,k}\simeq\iM(\Gamma_0(N))_{\overline{\chi},k}\]
where $\big(\VS{n\in\M}a_nq^{\frac nt}\big)^c=\VS{n\in\M}\overline{a_n}q^{\frac nt}$. Next is an expansion of Lemma \ref{V1}.
\begin{Lem}\label{V2}For $h\in\N$, we have
\[\iM(\Gamma_0(N))_{\chi,k}\subset\iM(\Gamma_0(hN))_{\chi,k},\]
\[\iM(\Gamma_0(N))_{\chi,k}^{\ang h}=\iM(\Gamma_0(hN))_{\chi,k}\cap\C[[q^h]].\]
\end{Lem}

\vspace{0.5cm}

Let ${\tt 1}_N\neq\chi\in{\rm Hom}(\Z/N^\times,\C^\times)$. We denote the first generalized Bernoulli number by $B_{\chi}$. If $\chi$ is primitive then $B_{\chi}=\frac1N\VS{a=1}^N\chi(a)a$.
Note
\[\chi(-1)\neq-1\Longrightarrow B_{\chi}=0.\]

When $\chi(-1)=-1$, put
\[\iE_{\chi}=1-\z\frac2{B_{\chi}}\VS{n\in\N}\VS{d|n}\chi(d)q^n\]
then $\iE_{\chi}\in\iM(\Gamma_0(N))_{\chi,1}$. For example
\begin{align*}
\iE_{\rho_4}&=1+4\VS{n\in\N}\VS{d|n}\rho_4(d)q^n\\
\iE_{\rho_4\rho_8}&=1+2\VS{n\in\N}\VS{d|n}\rho_4\rho_8(d)q^n\\
\iE_{\rho_4\chi_{16}}&=1+(1-i)\VS{n\in\N}\VS{d|n}\rho_4\chi_{16}(d)q^n\\
\iE_{\rho_4\chi_{32}}&=1+(i-\zeta_8^3)\VS{n\in\N}\VS{d|n}\rho_4\chi_{32}(d)q^n\\
\iE_{\rho_4\chi_{32}^5}&=1+(i+\zeta_8^3)\VS{n\in\N}\VS{d|n}\rho_4\chi_{32}^5(d)q^n
\end{align*}

Also put
\[\iG_{\chi}=\VS{n\in\N}\VS{d|n}\chi(d)q^n.\]
If $N=p^n$ ($p:\text{prime}$, $n\in\N$), then $\iG_{\chi}=\frac{B_{\chi}}2\big(\iE_{\chi}-\iE_{\chi}^{\ang p}\big)\in\iM(\Gamma_0(p^{n+1}))_{\chi,1}$.

Moreover let ${\tt 1}_M\neq\psi\in{\rm Hom}(\Z/M^\times,\C^\times)$. When $\psi\chi(-1)=-1$, put
\[\iG_{\psi|\chi}=\VS{n\in\N}\Big(\VS{d|n}\psi(d)\chi(\z\frac nd)\Big)q^n\]
then $\iG_{\psi|\chi}\in\iM(\Gamma_0(MN))_{\psi\chi,1}$. Note $\iG_{\chi|\psi}=\iG_{\psi|\chi}$.

See \cite[\S4]{DS} or \cite[\S5.3]{St} for further details. Every space of modular forms is spanned by Eisenstein series and cusp forms (\cite[Theorem 4.5.2, 4.6.2, 4.8.1]{DS}). We see many cusp forms of weight 1 at \cite{L1}.

\newpage

\subsection{The case $N=4$}

\begin{Lem}\label{rel22}$\frac12(\eta^{\natural2}+\eta^{\flat2})=\eta^{\natural\ang22}$ and $\frac18(\eta^{\natural2}-\eta^{\flat2})=\eta^{\sharp\ang22}$.\end{Lem}
\begin{proof}We see
\[\z\frac12(\eta^{\natural4}+\eta^{\flat4})=\iE_{2,2}=\iE_{2,2}^{\ang{\frac12}\ang2}=(2\eta^{\natural4}-\eta^{\flat4})^{\ang2}\]
and $\frac14(\eta^{\natural4}+2\eta^{\flat\ang24}+\eta^{\flat4})=\eta^{\natural\ang24}$.
\[\z\frac18(\eta^{\natural2}-\eta^{\flat2})=\dfrac{\frac1{16}(\eta^{\natural4}-\eta^{\flat4})}{\frac12(\eta^{\natural2}+\eta^{\flat2})}=\dfrac{\eta^{\sharp4}}{\eta^{\natural\ang22}}=\eta^{\sharp\ang22}.\]
\end{proof}
We also see
\[\z\frac12(\eta^{\natural}+\eta^{\flat})=\sqrt{\frac14\big(\eta^{\natural2}+\eta^{\flat2}+2\eta^{\flat\ang2}\big)}=\eta^{\natural\ang4},\]
\[\z\frac14(\eta^{\natural}-\eta^{\flat})=\dfrac{\frac18(\eta^{\natural2}-\eta^{\flat2})}{\frac12(\eta^{\natural}+\eta^{\flat})}=\dfrac{\eta^{\sharp\ang22}}{\eta^{\natural\ang4}}=\eta^{\sharp\ang4}.\]

Moreover $\iE_{2,2}=2\iE_{1,2}^{\ang2}-\iE_{1,2}$ and
\[\iE_{4,2}^{\ang{\frac12}}=\z\frac13\big(4\iE_{1,2}^{\ang2}-\iE_{1,2}^{\ang{\frac12}}\big)=\frac23\iE_{2,2}+\frac13\iE_{2,2}^{\ang{\frac12}}=\eta^{\natural4}.\]\\

Note $\big\{\iE_{4,2},\iE_{2,2}^{\ang2}\big\}$ is a basis of $\cM(\Gamma_0(4))_2=\cM(2)_2^{\ang2}$. We get $\iE_{4,2}=\iE_{\rho_4}^2$. SAGE command for this identity :

\verb~R.<q>=PowerSeriesRing(QQ);prec=20~

\verb~E12=eisenstein_series_qexp(2,prec,normalization='constant')~

\verb~E42=(4*E12(q^4)-E12)/3~

\verb~x=DirichletGroup(4).0~

\verb~def dsx(n):~

\verb~    return sum([x(t) for t in divisors(n)])~

\verb~Ex=1+4*sum([dsx(n)*q^n for n in [1..prec]])~

\verb~E42-Ex^2~

$\longrightarrow$ \verb~O(q^20)~\\

In particular $\iE_{\rho_4}^{\ang{\frac12}}=\eta^{\natural2}$. Since $\iE_{\rho_4}\in\iM(\Gamma_0(4))_{\rho_4,1}$ we get $\eta^{\natural2}\in\iM(2)_{\rho_4,1}$.

Since $\eta^{\natural2}\eta^{\flat2}=\eta^{\flat\ang24}\in\iM(\Gamma_0(2))_{c_4,2}$ we get $\eta^{\flat2}\in\iM(2)_{c_4\rho_4,1}$, and $\eta^{\sharp2}\in\iM(2)_{b_4\rho_4,1}$ as well where
\[b_4:(\BM a&b\\c&d\EM)\mapsto(-1)^{\frac b2},\6c_4:(\BM a&b\\c&d\EM)\mapsto(-1)^{\frac c2}.\]

Note $\eta^{\sharp2}=\frac14\sqrt{\eta^{\natural4}-\eta^{\flat4}}$ by Lemma \ref{rel21}.

\begin{Thm}$\natural:\C\big[\eta^{\natural2},\eta^{\flat2},\eta^{\sharp2}\big]^{[1,1,1]}\bto\cM(4)$ and $\cM(4)\eta^6=\cS(4)$.\end{Thm}
\begin{proof}Considering $\eta^{\sharp2}\big|(\BM1&2\\0&1\\\EM)=-\eta^{\sharp2}$, we get $\hto$ of the first assertion.

The dimension formula says for $k\in\M$.
\[\dim\cM(4)_k=2k+1,\6\dim\cS(4)_k=(\dim\cM(4)_k-6)^+=\dim\cM(4)_{k-3}.\]
\end{proof}

\newpage

\subsection{Half-integer weight modular form}

For $\kappa\in\frac12\M$ and a congruence group $\Gamma\subset{\rm SL}_2(\Z)$, we define
\[\cM(\Gamma)_\kappa=\big\{f\in\cO(\cH) \,\big|\, \forall\gamma\in\Gamma.f|_\kappa\gamma=f,\,\forall\alpha\in{\rm SL}_2(\Z).f|_\kappa\alpha:{\rm holomorphic\;at\;}\infty\big\}.\]
This natural definition may be different from usual one, but we can make the $\frac12\M$-graded ring
\[\cM(\Gamma)_{\frac12\M}=\VT{\kappa\in\frac12\M}\cM(\Gamma)_\kappa.\]
Note $\cM(\Gamma)=\cM(\Gamma)_{\frac12\M}\big|_\M$. Lemma \ref{cuspform} has continuation.
\begin{Lem}$\eta^3\in\cS(8)_{\frac32}$ and $\eta\in\cS(24)_{\frac12}$.\end{Lem}

Generally, for a set $X$ denote the group of all maps $X\to\C^\times$ by ${\rm Map}(X,\C^\times)$. We regard ${\rm Map}(\Z/N^\times,\C^\times)\subset{\rm Map}(\Z/hN^\times,\C^\times)$ naturally. For a congruence subgroup $\Gamma$ and for $\chi\in{\rm Map}(\Z/N^\times,\C^\times)$ we regard $\chi=\chi\circ d_N\in{\rm Map}(\Gamma,\C^\times)$.

For $\chi\in{\rm Map}(\Gamma,\C^\times)$, we define
\[\iM(\Gamma)_{\chi,\kappa}=\big\{f\in\cM(\ker\chi)_\kappa \,\big|\, f|_\kappa\gamma=\chi(\gamma)f\text{ for all }\gamma\in\Gamma\big\}.\]
Note
\[(\BM -1&0\\0&-1\EM)\in\Gamma,\;\chi(\BM -1&0\\0&-1\EM)\neq i^{-2\kappa} \Longrightarrow\iM(\Gamma)_{\chi,\kappa}=\{0\}.\]
For a subgroup $X\subset{\rm Map}(\Z/N,\C^\times)$ we may make the graded ring $\iM(\Gamma)_X$.
\\

Let $\sqrt{\rho_4}\in{\rm Map}(\Z/4^\times,\C^\times)$ by $1\mapsto1$ and $-1\mapsto i$.

\begin{Prop}\label{dec2}For a subgroup $G\subset\Z/N^\times$, if $\rho_4(G)=\{1\}$ then $\Gamma_G(N)\subset\Gamma_1(4)$ and
\[\cM(\Gamma_G(N))_{\frac12\M}=\iM(\Gamma_0(N))_{\ang{\sqrt{\rho_4}}H}\]
where $H=\{f\in{\rm Hom}(\Z/N^\times,\C^\times) \,|\, f(G)=\{1\}\}$.
\end{Prop}

For example, $\cM(\Gamma_{\ang5}(2^n))_{\frac k2}=\iM(\Gamma_0(2^n))_{\overline{\sqrt{\rho_4}}^k,\frac k2}$ for $n\in\M+2$.

Note
\[f\mapsto f^c:\iM(\Gamma_0(N))_{\overline{\sqrt{\rho_4}}^{2\kappa}\chi,\kappa}\simeq\iM(\Gamma_0(N))_{\overline{\sqrt{\rho_4}}^{2\kappa}\overline{\chi},\kappa}.\]\\

We don't construct $|\alpha:\cM(\Gamma)_{\frac12\M}\bto\cM(\Gamma\TL\alpha)_{\frac12\M}$ generally, however for $h\in\N$ and $\kappa\in\frac12\M$, if $h|c$ then
\[f^{\ang h}|_\kappa(\BM a&b\\c&d\EM)=(\z\frac hd)^{2\kappa}(f|_\kappa(\BM a&bh\\c/h&d\EM))^{\ang h}\]
since
\begin{align*}f^{\ang h}|_\kappa(\BM a&b\\c&d\EM)(z)&=((\z\frac cd)(cz+d)^{1/2})^{-2\kappa}f((\BM h&0\\0&1\EM)(\BM a&b\\c&d\EM)z)\\
&=(\z\frac hd)^{-2\kappa}((\frac{c/h}d)(cz+d)^{1/2})^{-2\kappa}f((\BM a&bh\\c/h&d\EM)(\BM h&0\\0&1\EM)z)\\
&=(\z\frac hd)^{-2\kappa}f|_\kappa(\BM a&bh\\c/h&d\EM)(hz).\end{align*}
Similarly if $h|b$ then
\[f^{\ang{\frac1h}}|_\kappa(\BM a&b\\c&d\EM)=(\z\frac hd)^{2\kappa}(f|_\kappa(\BM a&b/h\\ch&d\EM))^{\ang h}.\]

\newpage

\subsection{Theta function}

Let $\theta=\VS{n\in\Z}q^{n^2}$, then $\theta\in\iM(\Gamma_0(4))_{\overline{\sqrt{\rho_4}},\frac12}$. The dimension formula says $\dim\iM(\Gamma_0(4))_{\rho_4,1}=1$ and $\theta^2=\iE_{\rho_4}$, $\eta^{\natural}=\theta^{\ang{\frac12}}$. We also have
\[\eta^{\flat}=\eta^{\natural}|(\BM1&1\\0&1\\\EM)=\big(\VS{n\in\Z}(-1)^nq^{n^2}\big)^{\ang{\frac12}}.\]
These are well known identities due to Jacobi. Expanding out the relation $\theta^2=\iE_{\rho_4}$ gives Jacobi's two-square theorem : for $n\in\N$
\[\#\big\{(a,b)\in\Z\times\Z \,\big|\, a^2+b^2=n\big\}=4\VS{d|n}\rho_4(d).\]

Let $\xi_8=\sqrt{\rho_4}\rho_8$, then $\xi_8^2=\rho_4$ and
\[\eta^{\natural}\in\iM(2)_{\xi_8,\frac12},\]
\[\eta^{\flat}\in\iM(\Gamma_0(4)\cap\Gamma^0(2))_{c_8\xi_8,\frac12},\]
\[\eta^{\sharp}\in\iM(\Gamma_0(2)\cap\Gamma^0(4))_{b_8\xi_8,\frac12}\]
where
\[b_8:(\BM a&b\\c&d\EM)\mapsto(-1)^{\frac b4},\6c_8:(\BM a&b\\c&d\EM)\mapsto(-1)^{\frac c4}.\]
\\

For $\chi\in{\rm Hom}(\Z/N^\times,\C^\times)$ such that $\chi(-1)=1$, put
\[\theta_{\chi}=\VS{n\in\N}\chi(n)q^{n^2}\]
then $\theta_{\chi}\in\iM(\Gamma_0(4N^2))_{\overline{\sqrt{\rho_4}}\chi,\frac12}$. In particular $\eta=\theta_{\rho_4\rho_3}^{\ang{\frac1{24}}}$.

Remark $\theta_{{\tt 1}_p}=\frac12(\theta-\theta^{\ang{p^2}})$ for a prime $p$, in particular $\eta^{\sharp}=\frac12(\eta^{\natural}-\eta^{\natural\ang4})^{\ang{\frac14}}=\theta_{{\tt 1}_2}^{\ang{\frac18}}$.

In their paper, Serre and Stark prove that
$\cM(\Gamma_1(N))_{\frac12}$ is spanned by
\[\theta^{\ang{h}} \text{ with } 4h|N,\6 \theta_{\chi}^{\ang h} \text{ with } 4c_\chi^2h|N\]
where $c_\chi$ is the conductor of $\chi$. If $\chi$ is totally-even, that is, those $\chi$ whose prime-power components $\chi_p$ are all even, or $\chi$ is square of some character of conductor $\gcd(c_\chi,2)c_\chi$, then $\theta_{\chi}$ is (closely related to) a twist of $\theta$, and hence one would not expect it to be cuspidal. However, if $\chi$ is even but not totally even, then $\theta_{\chi}$ turns out to be a cusp form.

As for $\chi\in{\rm Hom}(\Z/N^\times,\C^\times)$ such that $\chi(-1)=-1$, put
\[\Theta_{\chi}=\VS{n\in\N}\chi(n)nq^{n^2}\]
then $\Theta_{\chi}\in\iS(\Gamma_0(4N^2))_{\overline{\sqrt{\rho_4}}\chi,\frac32}$ (proposition 2.2 of Shimura's Annals paper).

\newpage

\subsection{The case $N=8$}

Note $\eta^{\natural\ang{\frac14}}=\eta^{\natural}+2\eta^{\sharp}$, $\eta^{\flat\ang{\frac14}}=\eta^{\natural}-2\eta^{\sharp}$ and put
\[\eta^{\uparrow}=\eta^{\natural\ang{\frac12}}|(\BM 1&1\\0&1\EM)=\eta^{\natural\ang2}+2i\eta^{\sharp\ang2}\]
\[\eta^{\downarrow}=\eta^{\natural\ang{\frac12}}|(\BM 1&-1\\0&1\EM)=\eta^{\natural\ang2}-2i\eta^{\sharp\ang2}\]
then $\eta^{\uparrow},\eta^{\downarrow}\in\iM(4)_{\xi_8\rho_8,\frac12}$. Lemma \ref{rel22} leads to
\[\sqrt{\eta^{\uparrow}\eta^{\downarrow}}=\sqrt{\eta^{\natural\ang22}+4\eta^{\sharp\ang22}}=\z\sqrt{\frac12(\eta^{\natural2}+\eta^{\flat2})+4\cdot\frac18(\eta^{\natural2}-\eta^{\flat2})}=\eta^{\natural}.\]

Also note $\eta^{\natural\ang{\frac12}}=\frac{1-i}2\eta^{\uparrow}+\frac{1+i}2\eta^{\downarrow}$, $\eta^{\flat\ang{\frac12}}=\frac{1+i}2\eta^{\uparrow}+\frac{1-i}2\eta^{\downarrow}$ and
\[\eta^{\flat}=\sqrt{\eta^{\natural\ang{\frac12}}\eta^{\flat\ang{\frac12}}}=\z\sqrt{\frac12(\eta^{\uparrow2}+\eta^{\downarrow2})},\]
\[\eta^{\sharp}=\z\sqrt{\frac18(\eta^{\natural\ang{\frac12}2}-\eta^{\flat\ang{\frac12}2})}=\sqrt{\frac1{8i}(\eta^{\uparrow2}-\eta^{\downarrow2})}.\]

\begin{Thm}$\natural:\C\Big[\sqrt{\eta^{\uparrow}},\sqrt{\eta^{\downarrow}},\eta^{\flat},\eta^{\sharp}\Big]^{[\frac14,\frac14,\frac12,\frac12]}\Big|_{\frac12\M}\bto\cM(8)_{\frac12\M}=\iM(8)_{\ang{\xi_8,\rho_8},\frac12\M}$ and $\cM(8)_{\frac12\M}\eta^3=\cS(8)_{\frac12\M}$.
\end{Thm}
\begin{proof}Let $R=\C\Big[\sqrt{\eta^{\uparrow}},\sqrt{\eta^{\downarrow}}\Big]^{[\frac14,\frac14]}$ then $\natural:R|_{\frac12\M}\hto\iM(4)_{\ang{\xi_8,\rho_8},\frac12\M}$.

Considering $\eta^{\flat}\big|(\BM1&0\\4&1\\\EM)=-\eta^{\flat}$ and $\eta^{\sharp}\big|(\BM1&4\\0&1\\\EM)=-\eta^{\sharp}$, we get $\hto$ of the former assertion.

The dimension formula says
\[\dim\cM(8)_k=16k-4,\6\dim\cS(8)_k=\z\delta_{k,2}+16(k-\frac32)-4\]
for $k\in\N$. Note $\eta^3\in\cS(8)_{\frac32}$ and as for $k\in2\M+1$
\[\dim\cM(8)_{\frac k2}\leq\dim\cS(8)_{\frac{k+3}2}\leq\delta_{k,1}+8k-4.\]

We calculate
\[{\rm Dim\,LHS}|_{\frac12\M}=\frac{(1+t^{\frac12})^2}{(1-t^{\frac14})^2}\Big|_{\frac12\M}={\rm Dim}\,\cM(8)_{\frac12\M}.\]

On the latter assertion, we easily get $\subset$ and $\cM(8)_{\frac k2}\eta^3=\cS(8)_{\frac{k+3}2}$ for $k\in2\M+1$. As for $k\in2\N$
\[\cS(8)_{\frac{k+3}2}\subset\big\{f\in\cM(8)_{\frac{k+3}2} \,\big|\, f\eta^\flat,f\eta^\sharp\in\cS(8)_{\frac{k+4}2}\big\}\subset\big(\cM(8)_{\frac{k+1}2}\eta^\sharp\cap\cM(8)_{\frac{k+1}2}\eta^\flat\big)\eta^\natural.\]
Let
\[f\in\cM(8)_{\frac{k+1}2}\eta^\sharp\cap\cM(8)_{\frac{k+1}2}\eta^\flat\]
then there exist $x_0\in R_{\frac{k+1}2}$, $x_\sharp,y_\flat\in R_{\frac k2}$, $y_1\in R_{\frac{k-1}2}$ such that
\[f\in\big(x_0+x_\sharp\eta^\sharp+S_{\frac k2}\eta^\flat\big)\eta^\sharp\cap\big(S_{\frac{k+1}2}+y_\flat\eta^\flat+y_1\eta^\flat\eta^\sharp\big)\eta^\flat.\]
where $S_{\frac k2}=R_{\frac k2}+R_{\frac{k-1}2}\eta^\sharp$. We see $x_0=y_1\eta^{\flat2}$ and $x_\sharp\eta^{\sharp2}=y_\flat\eta^{\flat2}$ thus
\[f\in\cM(8)_{\frac k2}\eta^\sharp\eta^\flat.\]\end{proof}

The above Theorem derives $\natural:\C\big[\eta^{\natural},\eta^{\flat},\eta^{\sharp}\big]^{[\frac12,\frac12,\frac12]}\bto\iM(8)_{\ang{\xi_8},\frac12\M}$.

We have $\dim\cS(8)_{\frac32}=1$ and $\eta^3=\Theta_{\rho_4}^{\ang{\frac18}}$.

Also $\dim\iM(\Gamma_0(2)\cap\Gamma^0(4))_{\rho_4\rho_8,1}=1$ and $\eta^{\natural\ang{\frac12}}\eta^{\natural}=\iE_{\rho_4\rho_8}^{\ang{\frac14}}$.

\newpage

\section{For the case $N=16,32$}

\subsection{Rational weight}

When $f$\ is a modular form of weight 1, $\sqrt[n]f$ should be a form of weight $\frac 1n$. To treat it neatly, on a congruence group $\Gamma$, with $J:\Gamma\times\cH\to\cH$ satisfing
\[J((\BM a&b\\c&d\EM),z)^n=\z\frac1{cz+d}\]
we define
\[f|_{\frac kn,J}\alpha(z)=J(\alpha,z)^kf(\alpha z).\]
We easily see
\[f|_{\frac kn,J}\alpha\cdot f'|_{\frac{k'}n,J}\alpha=(ff')|_{\frac{k+k'}n,J}\alpha,\]
thus we can make the $\frac1n\M$-graded ring
\[\cM(\Gamma)_{\frac1n\M,J}=\VT{\kappa\in\frac1n\M}\cM(\Gamma)_{\kappa,J}.\]

For $\phi\in{\rm Map}(\Gamma,\C^\times)$ and $x\in\cO(\cH)$, if $x$ is non-zero on $\cH$ then define $J[\phi,x]:\Gamma\times\cH\to\cH$ by
\[J[\phi,x](\alpha,z)=\frac{\phi(\alpha)x(z)}{x(\alpha z)}.\]\\

Let $J[4]=J\big[\sqrt{\xi_8},\sqrt{\eta^{\natural}}\big]$ on $\Gamma(2)$ and $|_{\frac k4}=|_{\frac k4,J[4]}$ unless otherwise noted, then
\[\sqrt{\eta^{\natural}}\in\iM(2)_{\sqrt{\xi_8},\frac14}\]
\[\sqrt{\eta^{\flat}}\in\iM(\Gamma_0(8)\cap\Gamma^0(2))_{c_{16}\sqrt{\xi_8}\rho_8,\frac14}\]
\[\sqrt{\eta^{\sharp}}\in\iM(\Gamma_0(2)\cap\Gamma^0(8))_{b_{16}\sqrt{\xi_8}\rho_8,\frac14}\]
where
\[b_{16}:(\BM a&b\\c&d\EM)\mapsto(-1)^{\frac b8},\6c_{16}:(\BM a&b\\c&d\EM)\mapsto(-1)^{\frac c8}.\]

A basis of $\iM(4)_{\rho_8,2}$ is $\{\eta^{\uparrow},\eta^{\downarrow}\}\cdot\{\eta^{\uparrow},\eta^{\downarrow}\}\cdot\eta^{\natural}\cdot\{\eta^{\uparrow},\eta^{\downarrow}\}$ and $\sqrt{\eta^{\downarrow}\eta^{\natural\ang{\frac12}}\eta^{\natural}\eta^{\natural\ang2}}=\iE_{\rho_4\chi_{16}}^{\ang{\frac14}}$.
\begin{Lem}$\sqrt{\eta^{\downarrow}\eta^{\natural\ang{\frac12}}\eta^{\natural\ang2}}\in\iM(4)_{\sqrt{\xi}_8\chi_{16},\frac34}$.\end{Lem}

SAGE command for $\sqrt{\eta^{\downarrow\ang4}\eta^{\natural\ang2}\eta^{\natural\ang4}\eta^{\natural\ang8}}=\iE_{\rho_4\chi_{16}}$ :

\verb~K.<zeta4>=CyclotomicField(4);R.<q>=PowerSeriesRing(K);prec=20~

\verb~from sage.modular.etaproducts import qexp_eta~

\verb~E=qexp_eta(ZZ[['q']], prec);Ef=E^2/E(q^2);En=Ef(-q)~

\verb~x=(DirichletGroup(16).0)*(DirichletGroup(16).1)~

\verb~def dsx(n):~

\verb~    return sum([x(t) for t in divisors(n)])~

\verb~Ex=1+(1-zeta4)*sum([dsx(n)*q^n for n in [1..prec]])~

\verb~sqrt(En(zeta4^3*q)*En*En(q^2)*En(q^4))-Ex~

$\longrightarrow$ \verb~O(q^20)~

\newpage

\subsection{The case $N=16$}

\begin{Lem}Let $\big(\frac14\Z/\Z\times\frac18\M\big)$-graded ring
\[R=\C\Big[\sqrt[4]{\eta^{\uparrow}},\sqrt[4]{\eta^{\downarrow}},\sqrt{\eta^{\flat}},\sqrt{\eta^{\sharp}}\Big]^{[(\frac14,\frac18),(\frac34,\frac18),(0,\frac14),(0,\frac14)]}\]
and $X=\ang{\sqrt{\xi_8},\rho_8}$ then $\natural:R|_{\{0\}\times\frac14\M}\bto\iM(16)_{X,\frac14\M}$.\end{Lem}
\begin{proof}
First, we see
\[R|_{\{0\}\times\frac14\M}=\C\Big[\eta^{\uparrow},\eta^{\downarrow},\sqrt{\eta^{\natural}},\sqrt{\eta^{\flat}},\sqrt{\eta^{\sharp}}\Big]^{[\frac12,\frac12,\frac14,\frac14,\frac14]}\]
and so get $\hto$.

With $e^{\frac14}=t^{(\frac14,0)}$ and $\tau^{\frac14}=t^{(0,\frac14)}$, we calculate
\[{\rm Dim}\,R=\frac{(1-\tau)^2}{(1-e^{\frac14}\tau^{\frac18})(1-e^{\frac34}\tau^{\frac18})(1-\tau^{\frac14})^2}\]
\[{\rm Dim}\,R|_{\{0,\frac12\}\times\frac14\M}=\frac{(1+\tau^{\frac14})(1-\tau)^2}{(1-e^{\frac12}\tau^{\frac14})^2(1-\tau^{\frac14})^2}\]
\begin{align*}{\rm Dim}\,R|_{\{0\}\times\frac14\M}&=\frac{(1+\tau^{\frac14})(1+\tau^{\frac12})^3}{(1-\tau^{\frac14})^2}\\&=\frac{1+\tau^{\frac14}+3\tau^{\frac12}+3\tau^{\frac34}+3\tau+3\tau^{\frac54}+\tau^{\frac32}+\tau^{\frac74}}{(1-\tau^{\frac14})^2}\end{align*}

The dimension formula says $\dim\iM(16)_{X,k}=64k-40+3\delta_{1,k}$ for $k\in\N$.

For $k\in\M$, we see
\[R_{(0,\frac k4)}\eta^{\frac32}\subset\iM(16)_{X,\frac k4}\eta^{\frac32}\subset\iS(16)_{X,\frac{k+3}4}.\]

As for $k\in4\M+1$, the dimension formula says $\dim R_{(0,\frac k4)}=\dim\iS(16)_{X,\frac{k+3}4}$ hence above two $\subset$ are $=$.

For $k\in4\N$ and then for $k\in4\M+3$ and then for $k\in4\M+2$
\[\iS(16)_{X,\frac{k+3}4}\subset\big(R_{(0,\frac{k+1}4)}\sqrt{\eta^\sharp}\cap R_{(0,\frac{k+1}4)}\sqrt{\eta^\flat}\big)\sqrt{\eta^\natural}\subset R_{(0,\frac k4)}\eta^{\frac32}.\]
\end{proof}

We have $\dim\big\{f\in\iM(16)_{\xi_8,\frac12}\,\big|\,f|(\BM 1&8\\0&1\EM)=-f\big\}=1$ and $\sqrt{\eta^{\Diamond}}=\theta_{\rho_8}^{\ang{\frac1{16}}}$.

With $A=\{\eta^{\uparrow},\eta^{\downarrow},\eta^{\flat\ang2}\}\sqrt{\eta^{\flat}}\cup\{\eta^{\natural},\eta^{\flat},\eta^{\sharp}\}\sqrt{\eta^{\natural}}$,

$A\sqrt{\eta^{\natural}\eta^{\Diamond}}$ is a basis of $\big\{f\in\iS(16)_{\xi_8^3,\frac32}\,\big|\,f|(\BM 1&8\\0&1\EM)=-f\big\}$ and $\eta^{\natural2}\sqrt{\eta^{\Diamond}}=\Theta_{\rho_4\rho_8}^{\ang{\frac1{16}}}$.

\newpage

Note $\sqrt{\dfrac{\eta^{\flat\ang{\frac12}}}{\eta^{\natural\ang{\frac12}}}}=\dfrac{\eta^{\flat}}{\eta^{\natural\ang{\frac12}}}$, $\sqrt{\dfrac{\eta^{\sharp\ang2}}{\eta^{\natural\ang2}}}=\dfrac{\eta^{\sharp}}{\eta^{\natural\ang2}}$ and for $(\BM a&b\\c&d\EM)\in\Gamma_0(4)\cap\Gamma^0(4)$,
\[\sqrt{\frac{\eta^{\flat\ang{\frac12}}}{\eta^{\natural\ang{\frac12}}}}\Big|_0(\BM a&b\\c&d\EM)=(-1)^{\frac c4}\rho_8(d)\sqrt{\frac{\eta^{\flat\ang{\frac12}}}{\eta^{\natural\ang{\frac12}}}}\]
\[\sqrt{\frac{\eta^{\sharp\ang2}}{\eta^{\natural\ang2}}}\Big|_0(\BM a&b\\c&d\EM)=(-1)^{\frac b4}\rho_8(d)\sqrt{\frac{\eta^{\sharp\ang2}}{\eta^{\natural\ang2}}}\]

\begin{Thm}Let $\big(\frac14\Z/\Z\times\frac14\Z/\Z\times\frac18\M\big)$-graded ring
\[R=\C\Big[\sqrt[4]{\eta^{\uparrow}},\sqrt[4]{\eta^{\downarrow}},\sqrt[4]{\eta^{\natural\ang{\frac12}}},\sqrt[4]{\eta^{\flat\ang{\frac12}}},\sqrt[4]{\eta^{\natural\ang2}},\sqrt[4]{\eta^{\sharp\ang2}}\Big]\]
\[\6^{[(\frac14,\frac34,\frac18),(\frac34,\frac14,\frac18),(\frac14,0,\frac18),(\frac34,0,\frac18),(0,\frac14,\frac18),(0,\frac34,\frac18)]}\]
then $\natural:R|_{\{0\}\times\{0\}\times\frac14\M}\bto\cM(16)_{\frac14\M}=\iM(16)_{\ang{\xi_8,\chi_{16}},\frac14\M}$ and $\cM(16)_{\frac14\M}\eta^{\frac32}=\cS(16)_{\frac14\M}$.\\
\end{Thm}
\begin{proof}With $e^{\frac14}=t^{(\frac14,0,0)}$, $f^{\frac14}=t^{(0,\frac14,0)}$ and $\tau^{\frac18}=t^{(0,0,\frac18)}$ we calculate
\[{\rm Dim}\,R=\frac{(1-\tau^{\frac12})^4}{(1-e^{\frac14}f^{\frac34}\tau^{\frac18})(1-e^{\frac34}f^{\frac14}\tau^{\frac18})(1-e^{\frac14}\tau^{\frac18})(1-e^{\frac34}\tau^{\frac18})(1-f^{\frac14}\tau^{\frac18})(1-f^{\frac34}\tau^{\frac18})},\]
\[{\rm Dim}\,R|_{\ang{\frac14}\times\ang{\frac14}\times\frac14\M}=\frac{(1-\tau^{\frac12})^4}{(1-e^{\frac12}f^{\frac12}\tau^{\frac14})^2(1-e^{\frac12}\tau^{\frac14})^2(1-f^{\frac12}\tau^{\frac14})^2}\cdot\]
\[\6\big((1+\tau^{\frac14})(1+(EF+2)\tau^{\frac14}+\tau^{\frac12})+(e^{\frac14}f^{\frac34}+e^{\frac34}f^{\frac14})\tau^{\frac18}(E+F)(\tau^{\frac18}+\tau^{\frac38})\big)\]
where $E=e^{\frac14}+e^{\frac34}$ and $F=f^{\frac14}+f^{\frac34}$,
\begin{align*}{\rm Dim}\,R|_{\{0,\frac12\}\times\{0,\frac12\}\times\frac14\M}&=\frac{(1-\tau^{\frac12})^4(1+\tau^{\frac14})^3}{(1-e^{\frac12}f^{\frac12}\tau^{\frac14})^2(1-e^{\frac12}\tau^{\frac14})^2(1-f^{\frac12}\tau^{\frac14})^2}\\
&=\frac{(1+\Phi\tau^{\frac14}+\Phi\tau^{\frac12}+\tau^{\frac34})^2(1+\tau^{\frac14})}{(1-\tau^{\frac14})^2}\end{align*}
where $\Phi=e^{\frac12}+f^{\frac12}+e^{\frac12}f^{\frac12}$. Note $\Phi^2=3+2\Phi$ and
\begin{align*}{\rm Dim}\,R|_{\{0\}\times\{0\}\times\frac14\M}&=\frac{(1+3\tau^{\frac12}+8\tau^{\frac34}+3\tau+\tau^{\frac32})(1+\tau^{\frac14})}{(1-\tau^{\frac14})^2}\\
&=\frac{1+\tau^{\frac14}+3\tau^{\frac12}+11\tau^{\frac34}+11\tau+3\tau^{\frac54}+\tau^{\frac32}+\tau^{\frac74}}{(1-\tau^{\frac14})^2}.\end{align*}

The dimension formula says $\dim\cM(16)_{\frac k4}=32k-80+3\delta_{k,4}$ for $k\in4\N$.

As for $k\in4\N+1$, $\dim\cM(16)_{\frac k4}\leq\dim\cS(16)_{\frac{k+3}4}=32k-80+\delta_{k,5}$.

As for $k\in4\N+2$, we see
\begin{align*}\cM(16)_{\frac k4}\eta^{\frac32}&\subset\cS(16)_{\frac{k+3}4}\\
&\subset\big\{f\in R_{(0,0,\frac{k+3}4)}\,\big|\,f^2\in\iM(16)_{\ang{\sqrt{\xi_8},\rho_8},\frac k2}\eta^3\big\}\\
&\subset R_{(0,0,\frac k4)}\eta^{\frac32}.\end{align*}

At last, the above result holds also for $k\in4\N+3$.
\end{proof}

\newpage

\subsection{The case $N=32$}

We get for $(\BM a&b\\c&d\EM)\in\Gamma_0(16)$
\[\sqrt{\frac{\eta^{\flat\ang2}}{\eta^{\natural\ang2}}}\Big|_0(\BM a&b\\c&d\EM)=(-1)^{\frac c{16}}\rho_8(d)\sqrt{\frac{\eta^{\flat\ang2}}{\eta^{\natural\ang2}}}\]
Note $\sqrt{\dfrac{\eta^{\flat\ang{\frac12}}}{\eta^{\natural\ang{\frac12}}}}=\dfrac{\eta^{\flat}}{\eta^{\natural\ang{\frac12}}}$, $\sqrt{\dfrac{\eta^{\sharp\ang2}}{\eta^{\natural\ang2}}}=\dfrac{\eta^{\sharp}}{\eta^{\natural\ang2}}$ and for $(\BM a&b\\c&d\EM)\in\Gamma^0(16)$
\[\sqrt{\frac{\eta^{\sharp\ang{\frac12}}}{\eta^{\natural\ang{\frac12}}}}\Big|_0(\BM a&b\\c&d\EM)=(-1)^{\frac b{16}}\rho_8(d)\sqrt{\frac{\eta^{\sharp\ang{\frac12}}}{\eta^{\natural\ang{\frac12}}}}\]

\begin{Lem}Let $\big(\frac18\Z/\Z\times\frac18\Z/\Z\times\frac1{16}\M\big)$-graded ring
\[R=\C\Big[\sqrt[8]{\eta^{\uparrow}},\sqrt[8]{\eta^{\downarrow}},\sqrt[8]{\eta^{\natural\ang{\frac12}}},\sqrt[8]{\eta^{\flat\ang{\frac12}}},\sqrt[8]{\eta^{\natural\ang2}},\sqrt[8]{\eta^{\sharp\ang2}}\Big]\]
\[\6^{[(\frac18,\frac38,\frac1{16}),(\frac38,\frac18,\frac1{16}),(\frac18,0,\frac1{16}),(\frac38,0,\frac1{16}),(0,\frac18,\frac1{16}),(0,\frac38,\frac1{16})]}\]
and $X=\ang{\rho_4,\chi_{16}}$ then $\natural:R|_{\{0\}\times\{0\}\times\M}\bto\iM(32)_X$.\end{Lem}
\begin{proof}
With $e^{\frac18}=t^{(\frac18,0,0)}$, $f^{\frac18}=t^{(0,\frac18,0)}$ and $\tau^{\frac1{16}}=t^{(0,0,\frac1{16})}$ we calculate
\begin{align*}&{\rm Dim}\,R\\
&=\frac{(1-\tau^{\frac12})^4}{(1-e^{\frac18}f^{\frac38}\tau^{\frac1{16}})(1-e^{\frac38}f^{\frac18}\tau^{\frac1{16}})(1-e^{\frac18}\tau^{\frac1{16}})(1-e^{\frac38}\tau^{\frac1{16}})(1-f^{\frac18}\tau^{\frac1{16}})(1-f^{\frac38}\tau^{\frac1{16}})},\end{align*}
\begin{align*}&{\rm Dim}\,R|_{\ang{\frac14}\times\ang{\frac14}\times\frac18\M}\\
&=\frac{(1-\tau^{\frac12})^4(1+e^{\frac12}f^{\frac12}\tau^{\frac18})(1+f^{\frac12}\tau^{\frac18})(1+e^{\frac12}\tau^{\frac18})}{(1-e^{\frac14}f^{\frac34}\tau^{\frac18})(1-e^{\frac34}f^{\frac14}\tau^{\frac18})(1-e^{\frac14}\tau^{\frac18})(1-e^{\frac34}\tau^{\frac18})(1-f^{\frac14}\tau^{\frac18})(1-f^{\frac34}\tau^{\frac18})},\end{align*}
\begin{align*}&{\rm Dim}R|_{\{0,\frac12\}\times\{0,\frac12\}\times\frac14\M}\\
&=\frac{(1-\tau^{\frac12})^4(1+\Phi\tau^{\frac14})(1+\tau^{\frac14})^3}{(1-e^{\frac12}f^{\frac12}\tau^{\frac14})^2(1-e^{\frac12}\tau^{\frac14})^2(1-f^{\frac12}\tau^{\frac14})^2}\\
&=\frac{(1+\Phi\tau^{\frac14}+\Phi\tau^{\frac12}+\tau^{\frac34})^2(1+\Phi\tau^{\frac14})(1+\tau^{\frac14})}{(1-\tau^{\frac14})^2},\end{align*}
where $\Phi=e^{\frac12}+f^{\frac12}+e^{\frac12}f^{\frac12}$,
\begin{align*}&{\rm Dim}R|_{\{0\}\times\{0\}\times\frac14\M}\\
&=\frac{(1+9\tau^{\frac12}+20\tau^{\frac34}+15\tau+12\tau^{\frac54}+7\tau^{\frac32})(1+\tau^{\frac14})}{(1-\tau^{\frac14})^3}\\
&=\frac{1+\tau^{\frac14}+9\tau^{\frac12}+29\tau^{\frac34}+35\tau+27\tau^{\frac54}+19\tau^{\frac32}+7\tau^{\frac74}}{(1-\tau^{\frac14})^2}\end{align*}

The dimension formula says $\dim\iM(32)_{X,k}=512k-416$ for $k\in\M+2$.
\end{proof}

We have $\sqrt{\eta^{\downarrow}}\sqrt[4]{\eta^{\Diamond}}=\theta_{\chi_{16}}^{\ang{\frac1{32}}}$ and $\eta^{\uparrow2}\sqrt{\eta^{\downarrow}}\sqrt[4]{\eta^{\Diamond}}=\Theta_{\rho_4\chi_{16}}^{\ang{\frac1{32}}}$.

\newpage

We see
\[\sqrt{\big(\sqrt{\eta^{\uparrow\ang{\frac12}}\eta^{\flat\ang{\frac14}}}-2\zeta_8^3\sqrt{\eta^{\uparrow}\eta^{\sharp\ang2}}\big)\sqrt{\eta^{\uparrow\ang{\frac12}}\eta^{\uparrow}\eta^{\natural\ang{\frac14}}\eta^{\natural\ang{\frac12}}\eta^{\natural}\eta^{\natural\ang2}}}=\iE_{\rho_4\chi_{32}}^{\ang{\frac18}}\]
\[\sqrt{\big(\sqrt{\eta^{\uparrow\ang{\frac12}}\eta^{\flat\ang{\frac14}}}+2\zeta_8^3\sqrt{\eta^{\uparrow}\eta^{\sharp\ang2}}\big)\sqrt{\eta^{\uparrow\ang{\frac12}}\eta^{\uparrow}\eta^{\natural\ang{\frac14}}\eta^{\natural\ang{\frac12}}\eta^{\natural}\eta^{\natural\ang2}}}=\iE_{\rho_4\chi_{32}^5}^{\ang{\frac18}}\]\\

Remark
\begin{align*}&\z\frac12\big(\sqrt{\eta^{\uparrow\ang{\frac12}}\eta^{\natural\ang{\frac14}}}+\sqrt{\eta^{\downarrow\ang{\frac12}}\eta^{\flat\ang{\frac14}}}\big)\\
&=\z\frac12\sqrt{(\eta^{\natural}+2i\eta^{\sharp})(\eta^{\natural}+2\eta^{\sharp})+2\eta^{\natural\ang{\frac12}}\eta^{\flat\ang{\frac12}}+(\eta^{\natural}-2i\eta^{\sharp})(\eta^{\natural}-2\eta^{\sharp})}\\
&=\z\sqrt{\frac12(\eta^{\natural2}+\eta^{\flat2})+2i\eta^{\sharp2}}\\
&=\sqrt{\eta^{\uparrow}\eta^{\natural\ang2}},\end{align*}
\begin{align*}&\z\frac{1-i}4\big(\sqrt{\eta^{\uparrow\ang{\frac12}}\eta^{\natural\ang{\frac14}}}-\sqrt{\eta^{\downarrow\ang{\frac12}}\eta^{\flat\ang{\frac14}}}\big)\\
&=\z\frac{1-i}4\sqrt{(\eta^{\natural}+2i\eta^{\sharp})(\eta^{\natural}+2\eta^{\sharp})-2\eta^{\natural\ang{\frac12}}\eta^{\flat\ang{\frac12}}+(\eta^{\natural}-2i\eta^{\sharp})(\eta^{\natural}-2\eta^{\sharp})}\\
&=\z\sqrt{\frac{-i}4(\eta^{\natural2}-\eta^{\flat2})+\eta^{\sharp2}}\\
&=\sqrt{\eta^{\downarrow}\eta^{\sharp\ang2}}.\end{align*}

Remark
\begin{align*}&\z\frac{1+i}2\sqrt{\eta^{\uparrow\ang{\frac12}}\eta^{\natural\ang{\frac14}}}+\frac{1-i}2\sqrt{\eta^{\downarrow\ang{\frac12}}\eta^{\flat\ang{\frac14}}}\\
&=\z\sqrt{\frac i2(\eta^{\natural}+2i\eta^{\sharp})(\eta^{\natural}+2\eta^{\sharp})+\frac12\eta^{\natural\ang{\frac12}}\eta^{\flat\ang{\frac12}}-\frac i2(\eta^{\natural}-2i\eta^{\sharp})(\eta^{\natural}-2\eta^{\sharp})}\\
&=\sqrt{\eta^{\flat2}-2(1-i)\eta^{\natural}\eta^{\sharp}},\end{align*}
\[\z\frac{1-i}2\sqrt{\eta^{\uparrow\ang{\frac12}}\eta^{\natural\ang{\frac14}}}+\frac{1+i}2\sqrt{\eta^{\downarrow\ang{\frac12}}\eta^{\flat\ang{\frac14}}}=\sqrt{\eta^{\flat2}+2(1-i)\eta^{\natural}\eta^{\sharp}}\]
\[\z\sqrt{\eta^{\flat2}-2(1-i)\eta^{\natural}\eta^{\sharp}}\sqrt{\eta^{\flat2}+2(1-i)\eta^{\natural}\eta^{\sharp}}=\eta^{\natural2}+4i\eta^{\sharp2}.\]

\newpage

\section{For other $N$}

\subsection{Preparations}

For a prime $p\neq2$ we define $\rho_p=(\z\frac\bullet p)$. Note $\rho_3=(\frac{-3}\bullet)$ and\[{\rm Hom}(\Z/3^\times,\C^\times)=\ang{\rho_3}.\]

Also define $\chi_N\in{\rm Hom}(\Z/N^\times,\C^\times)$ for $N=p^r$ with prime $p\neq2$ and $N\neq3$ by $\chi_N(g)=\zeta_{p^r-p^{r-1}}$ where $g$ is the minimum generator of $\Z/N^\times\simeq\Z/(p^r-p^{r-1})$. For example, $\chi_5(2)=i$, $\chi_7(3)=\zeta_6$ and $\chi_9(2)=\zeta_6$. Note $\chi_5^2=\rho_5$, $\chi_7^3=\rho_7$ and $\chi_9^3=\rho_3$. We see
\[{\rm Hom}(\Z/N^\times,\C^\times)=\ang{\chi_N}.\]
In addition
\[N\in2\N+1 \Longrightarrow {\rm Hom}(\Z/2N^\times,\C^\times)={\rm Hom}(\Z/N^\times,\C^\times)\]\[{\rm Hom}(\Z/12^\times,\C^\times)=\ang{\rho_4,\rho_3}\]\\

On some Eisenstein series
\begin{align*}
\iE_{\rho_3}&=1+6\VS{n\in\N}\VS{d|n}\rho_3(d)q^n\\
\iE_{\chi_5}&=1+(3-i)\VS{n\in\N}\VS{d|n}\chi_5(d)q^n\\
\iE_{\rho_7}&=1+2\VS{n\in\N}\VS{d|n}\rho_7(d)q^n\\
\iE_{\chi_7}&=1+(1-2\omega)\VS{n\in\N}\VS{d|n}\chi_7(d)q^n\\
\iE_{\chi_9}&=1+(1-\omega)\VS{n\in\N}\VS{d|n}\chi_9(d)q^n\\
\end{align*}

The dimension formula says $\dim\cM(\Gamma_0(p))_2=1$ for $p=3,5,7,13$ thus
\[\iE_{3,2}=\iE_{\rho_3}^2,\]
\[\iE_{5,2}=\iE_{\chi_5}\iE_{\overline{\chi_5}},\]
\[\iE_{7,2}=\iE_{\rho_7}^2=\iE_{\chi_7}\iE_{\overline{\chi_7}},\]
\[\iE_{13,2}=\iE_{\chi_{13}}\iE_{\overline{\chi_{13}}}=\iE_{\chi_{13}^3}\iE_{\overline{\chi_{13}^3}}=\iE_{\chi_{13}^5}\iE_{\overline{\chi_{13}^5}}.\]

\newpage

We abbreviate $\bot={/3/}$ and $\top="3"$. We see $f^{\bot}f^{\bot\top}=f^{\bot\ang33}$, $f^{\top}f^{\top\bot}=f^{\top\ang{\frac13}3}$ and $f^{\bot}f^{\top}f^{\bot\top}=f^8$.

\begin{Lem}\label{lev3}$\eta^{\bot}\in\iM(\Gamma^0(3))_{c_3\rho_3,1}$ and $\eta^{\top}\in\iM(\Gamma_0(3))_{b_3\rho_3,1}$ where
\[b_3:(\BM a&b\\c&d\EM)\mapsto\omega^{bd},\6c_3:(\BM a&b\\c&d\EM)\mapsto\overline\omega^{cd}.\]
\end{Lem}
\begin{proof}For $(\BM a&b\\c&d \EM)\in\Gamma^0(3)$, we see
\[\frac{(\frac3d)e^{\frac{3\pi i}{12}(3ac+\frac b3d-3cd-3acd^2+3d-3)}}{e^{\frac{\pi i}{12}(ac+bd-cd-acd^2+3d-3)}}=\z(\frac{-1}d)(\frac{-3}d)e^{\frac{2\pi i}3(a-d-ad^2)c+\frac{\pi i}2(d-1)}=\overline\omega^{cd}\rho_3(d),\]
\[\frac{(\frac3d)(\frac d{3c})(\frac{3c}d)e^{\frac{3\pi i}{12}(3ac+\frac b3d+3cd-3bc^2d-9c)}}{(\frac dc)(\frac cd)e^{\frac{\pi i}{12}(ac+bd+cd-bc^2d-3c)}}=\z(\frac d3)e^{\frac{2\pi i}3(a+d-bcd)c}=\overline\omega^{cd}\rho_3(d).\]

The latter assertion follows from the first one and $\eta^{\top}|_1(\BM 0&-1\\1&0\EM)=\frac1{\sqrt3^3i}\eta^{\bot}$.
\end{proof}

\vspace{0.5cm}

\begin{Lem}$\eta^{/5/}\in\iM(\Gamma^0(5))_{\rho_5,2}$ and $\eta^{"5"}\in\iM(\Gamma_0(5))_{\rho_5,2}$.\end{Lem}
\begin{proof}\[\frac{(\frac5d)e^{\frac{5\pi i}{12}(5ac+\frac b5d-5cd-5acd^2+3d-3)}}{e^{\frac{\pi i}{12}(ac+bd-cd-acd^2+3d-3)}}=\z(\frac5d)e^{\pi i(d-1)},\]
\[\frac{(\frac5d)(\frac d{5c})(\frac{5c}d)e^{\frac{5\pi i}{12}(5ac+\frac b5d+5cd-5bc^2d-15c)}}{(\frac dc)(\frac cd)e^{\frac{\pi i}{12}(ac+bd+cd-bc^2d-3c)}}=\z(\frac d5).\]

The latter assertion follows from the first one and $\eta^{"5"}|_2(\BM 0&-1\\1&0\EM)=-\frac1{\sqrt5^5}\eta^{/5/}$.\end{proof}

\vspace{0.5cm}

\begin{Lem}$\eta^{/7/}\in\iM(\Gamma^0(7))_{\rho_7,3}$ and $\eta^{"7"}\in\iM(\Gamma_0(7))_{\rho_7,3}$.\end{Lem}
\begin{proof}
\[\frac{(\frac7d)e^{\frac{7\pi i}{12}(7ac+\frac b7d-7cd-7acd^2+3d-3)}}{e^{\frac{\pi i}{12}(ac+bd-cd-acd^2+3d-3)}}=\z(\frac{-1}d)(\frac{-7}d)e^{\frac{\pi i}2(d-1)},\]
\[\frac{(\frac7d)(\frac d{7c})(\frac{7c}d)e^{\frac{7\pi i}{12}(7ac+\frac b7d+7cd-7bc^2d-21c)}}{(\frac dc)(\frac cd)e^{\frac{\pi i}{12}(ac+bd+cd-bc^2d-3c)}}=\z(\frac d7).\]

The latter assertion follows from the first one and $\eta^{"7"}|_3(\BM 0&-1\\1&0\EM)=\frac i{\sqrt7^7}\eta^{/7/}$.\end{proof}

\newpage

\subsection{The case $N=3$}

\begin{Thm}$\natural:\C\big[\eta^{\bot},\eta^{\top}\big]^{[1,1]}\bto\cM(3)$ and $\cM(3)\eta^8=\cS(3)$.
\end{Thm}
\begin{proof}The dimension formula says for $k\in\M$
\[\dim\cM(3)_k=k+1,\6\dim\cS(3)_k=(\dim\cM(3)_k-4)^+=\dim\cM(3)_{k-4}.\]
\end{proof}

\vspace{0.5cm}

Put $f^{\nwarrow}=f^{\bot}|(\BM 1&1\\0&1\EM)$ and $f^{\swarrow}=f^{\bot}|(\BM 1&-1\\0&1\EM)$ then
\begin{align*}\eta^{\nwarrow}\eta^{\swarrow}&=\VP{n\in\frac13\N\setminus\N}(1-\omega q^n)^3(1-\overline\omega q^n)^3\VP{n\in\N}(1-q^n)^4\\
&=\frac{\VP{n\in\N\setminus3\N}(1-q^n)^7\VP{n\in3\N}(1-q^n)^4}{\VP{n\in\frac13\N\setminus\N}(1-q^n)^3}=\eta^{\bot\top}\end{align*}
and $\eta^{\nwarrow}|_1(\BM 0&-1\\1&0\EM)=\overline\omega\eta^{\swarrow}$ since
\[(\BM 1&1\\0&1\EM)(\BM 0&-1\\1&0\EM)=(\BM 0&1\\-1&0\EM)(\BM 1&-1\\0&1\EM)(\BM 0&-1\\1&0\EM)(\BM 1&-1\\0&1\EM).\]

We see $\big\{f\in\cM(3)_1 \,\big|\, f|_1(\BM 1&1\\0&1\EM)=\overline\omega f\big\}=\{0\}$ and
\[\eta^{\bot}+\omega\eta^{\nwarrow}+\overline\omega\eta^{\swarrow}=0.\]
Acting $(\BM 0&-1\\1&0\EM)$ on the above identity, we get
\[\z3\eta^{\top}+\frac{1}{\sqrt3i}\eta^{\nwarrow}+\frac{-1}{\sqrt3i}\eta^{\swarrow}=0.\]\\

We have $\dim\cM(3)_1=2$ and $\eta^{\top}=\z\frac16\big(\iE_{\rho_3}^{\ang{\frac13}}-\iE_{\rho_3}\big)$.

We also see $\big\{f\in\cM(3)_1 \,\big|\, f|_1(\BM 1&1\\0&1\EM)=f\big\}=\C\iE_{\rho_3}$ and
\[\z\frac13(\eta^{\bot}+\eta^{\nwarrow}+\eta^{\swarrow})=\iE_{\rho_3}.\]
Moreover
\begin{align*}\z\sqrt3i\iE_{\rho_3}|_1(\BM 0&-1\\1&0\EM)&=\z\frac{-1}{\sqrt3i}(\eta^{\bot}+\eta^{\nwarrow}+\eta^{\swarrow})|_1(\BM 0&-1\\1&0\EM)\\
&=\z 3\eta^{\top}+\frac{\omega}{\sqrt3i}\eta^{\nwarrow}+\frac{-\overline\omega}{\sqrt3i}\eta^{\swarrow}\\
&=\z 6\eta^{\top}+\frac13(\eta^{\bot}+\eta^{\nwarrow}+\eta^{\swarrow})\\
&\6\z -\frac13(\eta^{\bot}+\omega\eta^{\nwarrow}+\overline\omega\eta^{\swarrow})-(3\eta^{\top}+\frac{1}{\sqrt3i}\eta^{\nwarrow}+\frac{-1}{\sqrt3i}\eta^{\swarrow})\\
&=\iE_{\rho_3}^{\ang{\frac13}}.\end{align*}

Note $\iE_{3,2}=\frac12(3\iE_{1,2}^{\ang3}-\iE_{1,2})$ and
\begin{align*}\iE_{9,2}^{\ang{\frac13}}&=\z\frac18(9\iE_{1,2}^{\ang3}-\iE_{1,2}^{\ang{\frac13}})=\frac34\iE_{3,2}-\frac14\iE_{3,2}^{\ang{\frac13}}\\
&=\z\frac34\big(\frac{1-\omega}3\eta^{\nwarrow}+\frac{1-\overline\omega}3\eta^{\swarrow}\big)^2-\frac14\big(\frac{\omega-1}{\sqrt3i}\eta^{\nwarrow}+\frac{-\overline\omega+1}{\sqrt3i}\eta^{\swarrow}\big)^2=\eta^{\nwarrow}\eta^{\swarrow}.\end{align*}

Put $\iE_{\rho_3}^{\nwarrow}=\iE_{\rho_3}^{\ang{\frac13}}|(\BM 1&1\\0&1\EM)$ and $\iE_{\rho_3}^{\swarrow}=\iE_{\rho_3}^{\ang{\frac13}}|(\BM 1&-1\\0&1\EM)$, then $\iE_{\rho_3}^{\nwarrow}|_1(\BM 0&-1\\1&0\EM)=\iE_{\rho_3}^{\swarrow}$ thus 
\[\iE_{1,4}=\iE_{\rho_3}\iE_{\rho_3}^{\ang{\frac13}}\iE_{\rho_3}^{\nwarrow}\iE_{\rho_3}^{\swarrow}.\]

\newpage

\subsection{The case $N=9$}

Note $\{\eta^{\nwarrow},\eta^{\top}\}\cdot\{\eta^{\nwarrow},\eta^{\top}\}\cdot\{\eta^{\swarrow},\eta^{\top}\}$ is a basis of $\cM(3)_3$ and $\sqrt[3]{\eta^{\nwarrow2}\eta^{\swarrow}}=\iE_{\chi_9}^{\ang{\frac13}}$.

Put
\[b_9:(\BM a&b\\c&d\EM)\mapsto\omega^{\frac{bd}3},\6c_9:(\BM a&b\\c&d\EM)\mapsto\overline\omega^{\frac{cd}3}.\]
Let $J[3]=J\big[c_9\rho_3,\sqrt[3]{\eta^\bot}\big]$ on $\Gamma_0(3)\cap\Gamma^0(3)$ and $|_{\frac k3}=|_{\frac k3,J[3]}$ unless otherwise noted. Since
\[\sqrt[3]{\eta^{\bot}}\in\iM(3)_{c_9\rho_3,\frac13}\]
and $\sqrt[3]{\eta^{\nwarrow}\eta^{\bot2}}=\dfrac{\eta^{\bot\ang32}}{\sqrt[3]{\eta^{\swarrow2}\eta^{\nwarrow}}}$ we see
\[\sqrt[3]{\eta^{\nwarrow}}\in\iM(3)_{\chi_9,\frac13}.\]
And
\[\sqrt[3]{\eta^{\swarrow}}\in\iM(3)_{\overline{\chi_9},\frac13},\]
\[\sqrt[3]{\eta^{\top}}\in\iM(3)_{b_9\rho_3,\frac13}.\]

\begin{Thm}Let
\[R=\C\Big[\sqrt[3]{\eta^{\nwarrow}},\sqrt[3]{\eta^{\swarrow}},\sqrt[3]{\eta^{\bot}},\sqrt[3]{\eta^{\top}}\Big]^{[\frac13,\frac13,\frac13,\frac13]}\]
then $\natural:R\bto\cM(9)_{\frac13\M}=\iM(9)_{\ang{\chi_9},\frac13\M}$. And $\cM(9)_{\frac13\M}\eta^{\frac83}=\cS(9)_{\frac13\M}$.\end{Thm}
\begin{proof}Considering $\sqrt[3]{\eta^{\bot}}\big|(\BM 1&0\\3&1 \EM)=\overline\omega\sqrt[3]{\eta^{\bot}}$ and $\sqrt[3]{\eta^{\top}}\big|(\BM 1&3\\0&1 \EM)=\omega\sqrt[3]{\eta^{\top}}$ we get $\hto$. Indeed
\[f_0+f_1\sqrt[3]{\eta^{\top}}+f_2\sqrt[3]{\eta^{\top}}^2=0 \Longrightarrow f_1\sqrt[3]{\eta^{\top}}=0+\overline\omega\cdot0|(\BM 1&3\\0&1 \EM)+\omega\cdot0|(\BM 1&6\\0&1 \EM)=0\]
for $f_0,f_1,f_2\in\C\Big[\sqrt[3]{\eta^{\nwarrow}},\sqrt[3]{\eta^{\swarrow}},\sqrt[3]{\eta^{\bot}}\Big]$.

We calculate
\[{\rm Dim}\,R=\dfrac{(1+t^{\frac13}+t^{\frac23})^2}{(1-t^{\frac13})^2}=\dfrac{1+2t^{\frac13}+3t^{\frac23}+2t+t^{\frac43}}{(1-t^{\frac13})^2}.\]

The dimension formuls says $\dim\cM(9)_k=27k-9$ for $k\in\N$.

For $k\in\M$, we see
\[R_{\frac k3}\eta^{\frac83}\subset\cM(9)_{\frac k3}\eta^{\frac83}\subset\cS(9)_{\frac{k+4}3}.\]

As for $k\in3\M+2$, the dimension formula says $\dim R_{\frac k3}=\dim\cS(9)_{\frac{k+4}3}$, hence above two $\subset$ are $=$. And then, for $k\in3\M+1$
\[\cS(9)_{\frac{k+4}3}\subset\big(R_{\frac{k+1}3}\sqrt[3]{\eta^{\top}}\cap R_{\frac{k+1}3}\sqrt[3]{\eta^{\bot}}\big)\sqrt[3]{\eta^{\nwarrow}\eta^{\swarrow}}\subset R\eta^{\frac83}.\]
\end{proof}

\newpage

\subsection{The case $N=27$}

Put $X=\eta^{\bot\ang3}+3\eta^{\top}$ then
\[\sqrt[3]{\eta^{\nwarrow}\eta^{\swarrow}}=\sqrt[3]{\eta^{\bot2}}+3\sqrt[3]{\eta^{\top\ang3}X},\]
\[\sqrt[3]{\eta^{\nwarrow\ang3}\eta^{\swarrow\ang3}}=\sqrt[3]{\eta^{\bot}X}+3\sqrt[3]{\eta^{\top\ang32}}.\]
\begin{Lem}Let
\[R=\C\Big[\sqrt[3]{\eta^{\bot}},\sqrt[3]{\eta^{\top}},\sqrt[3]{\eta^{\bot\ang3}},\sqrt[3]{\eta^{\top\ang3}},\sqrt[3]{X}\Big]^{[(0,\frac13),(0,\frac13),(\frac13,\frac13),(\frac13,\frac13),(\frac23,\frac13)]}\]
then $\natural:R_{\{0\}\times\M}\bto\iM(\Gamma_0(27)\cap\Gamma^0(9))_{\ang{\rho_3}}$.
\end{Lem}
\begin{proof}The dimension formula says $\dim\iM(\Gamma_0(27)\cap\Gamma^0(9))_{\ang{\rho_3},k}=27k-18+\delta_{k,1}$ for $k\in\N$. We calculate
\[{\rm Dim}\,R=\dfrac{(1-t)^3}{(1-t^{\frac13})^2(1-e^{\frac13}t^{\frac13})^2(1-e^{\frac23}t^{\frac23})},\]
\[{\rm Dim}\,R|_{\{0\}\times\frac13\M}=\dfrac{1+2t^{\frac23}+2t+3t^{\frac43}+t^{\frac53}}{(1-t^{\frac13})^2}.\]
\end{proof}

Note $\sqrt[3]{\eta^{\bot\ang3}\eta^{\top}}=\sqrt[3]{\eta\eta^{\ang3}}$ and
\[\iS(\Gamma_0(27)\cap\Gamma^0(9))_{\ang{\rho_3}}=R\sqrt[3]{\eta\eta^{\ang3}}\big|_{\{0\}\times\M}.\]
In particular $\sqrt[3]{\eta^{\bot\ang3}\eta^{\top}X}$ is the eigenform of $\iM(\Gamma_0(27)\cap\Gamma^0(9))_{\rho_3,1}$.

\newpage

\subsection{The case $N=6$}

Lemma \ref{lev3} easily derives
\[\eta^{\flat\bot}\in\iM(\Gamma^0(6))_{\rho_3,1},\6\eta^{\sharp\bot}\in\iM(\Gamma_0(2)\cap\Gamma^0(3))_{\rho_3,1}\]
\[\eta^{\flat\top}\in\iM(\Gamma_0(3)\cap\Gamma^0(2))_{\rho_3,1},\6\eta^{\sharp\top}\in\iM(\Gamma_0(6))_{\rho_3,1}\]
Indeed $\overline{\omega}^{\frac c2d}=(\overline{\omega}^{\frac c2d})^4=\overline{\omega}^{2cd}$ for $c\in2\Z$.\\[0.5cm]

We see $(\BM 2&0\\0&1\EM)(\BM 0&-1\\1&0\EM)=(\BM 0&-1\\1&0\EM)(\BM 1&0\\0&2\EM)$ and
\[\eta^{\sharp\top}|_1(\BM 0&-1\\1&0\EM)=\frac{\eta^{\top\ang22}|_2(\BM 0&-1\\1&0\EM)}{\eta^{\top}|_1(\BM 0&-1\\1&0\EM)}=\frac{(\frac1{2\sqrt3^3i}\eta^{\bot\ang{\frac12}})^2}{\frac1{\sqrt3^3i}\eta^{\bot}}=\frac1{4\sqrt3^3i}\eta^{\flat\bot}.\]
Note $\iE_{\rho_3}|_1(\BM 0&-1\\1&0\EM)=\frac1{\sqrt3i}\iE_{\rho_3}^{\ang{\frac13}}$ and $\iE_{\rho_3}^{\ang2}|_1(\BM 0&-1\\1&0\EM)=\frac1{2\sqrt3i}\iE_{\rho_3}^{\ang{\frac16}}$.

We have $\dim\iM(\Gamma_0(6))_{\rho_3,1}=2$ and
\[\z\eta^{\sharp\top}=\frac16(\iE_{\rho_3}-\iE_{\rho_3}^{\ang2}),\6\eta^{\flat\bot}=\big({-\iE_{\rho_3}}+2\iE_{\rho_3}^{\ang2}\big)^{\ang{\frac16}}.\]\\

Also
\[\eta^{\flat\top}|_1(\BM 0&-1\\1&0\EM)=\frac{\eta^{\top\ang{\frac12}2}|_2(\BM 0&-1\\1&0\EM)}{\eta^{\top}|_1(\BM 0&-1\\1&0\EM)}=\frac{(\frac2{\sqrt3^3i}\eta^{\bot\ang2})^2}{\frac1{\sqrt3^3i}\eta^{\bot}}=\frac4{\sqrt3^3i}\eta^{\sharp\bot}.\]

Take $x$ such that $\eta^{\flat\top}=x\iE_{\rho_3}^{\ang{\frac12}}+(1-x)\iE_{\rho_3}$, then $\eta^{\sharp\bot}=\frac32x\iE_{\rho_3}^{\ang{\frac23}}+\frac34(1-x)\iE_{\rho_3}^{\ang{\frac13}}$. Computing the coefficient of $q^0$ we get $x=\frac13$, therefore 
\[\z\eta^{\flat\top}=\frac13\big(\iE_{\rho_3}+2\iE_{\rho_3}^{\ang2}\big)^{\ang{\frac12}},\6\eta^{\sharp\bot}=\frac12\big(\iE_{\rho_3}+\iE_{\rho_3}^{\ang2}\big)^{\ang{\frac13}}.\]

\begin{Lem}We see
\[\sqrt{\eta^{\flat\top}}=\sqrt{\eta^{\sharp\bot\ang3}}+\sqrt{\eta^{\sharp\top}},\6\sqrt{\eta^{\natural\top}}=\sqrt{\eta^{\sharp\bot\ang3}}-\sqrt{\eta^{\sharp\top}},\]
\[\sqrt{\eta^{\flat\bot\ang3}}=\sqrt{\eta^{\sharp\bot\ang3}}+3\sqrt{\eta^{\sharp\top}},\6\sqrt{\eta^{\natural\bot\ang3}}=\sqrt{\eta^{\sharp\bot\ang3}}-3\sqrt{\eta^{\sharp\top}}.\]\end{Lem}
\begin{proof}First we see $\sqrt{\eta^{\sharp\bot\ang3}\eta^{\sharp\top}}=\eta^{\sharp}\eta^{\sharp\ang3}\in\cM(6)_1\cap\C[[q]]q^{\frac12}$.

The dimension formula says $\dim\cM(6)_1=6$ and we get $\eta^{\sharp}\eta^{\sharp\ang3}=\eta^{\sharp\top\ang{\frac12}}+\eta^{\sharp\top}$.

The first identity follows from $\eta^{\flat\top}=\big(\sqrt{\eta^{\sharp\bot\ang3}}+\sqrt{\eta^{\sharp\top}}\big)^2$. Acting $\eta^\natural=\eta^{\flat}|(\BM 1&1\\0&1\EM)$, the second one follows from $\sqrt{\eta^{\sharp\bot\ang3}}\in\C[[q]]$ and $\sqrt{\eta^{\sharp\top}}\in\C[[q]]q^{\frac12}$.
\end{proof}

Now, let $J[\natural]=J\big[\xi_8^3,\sqrt[3]{\eta^{\natural}}\big]$ on $\Gamma_0(2)\cap\Gamma^0(2)$, then $|_{\frac k2,J[\natural]}=|_{\frac k2}$. Since
\[\sqrt[3]{\eta^{\natural}}\in\iM(2)_{\xi_8^3,\frac16,J[\natural]}.\]
and $\sqrt[3]{\eta^{\natural\bot}\eta^{\natural}}=\eta^{\natural\ang{\frac13}}$, $\sqrt[3]{\eta^{\natural\top}\eta^{\natural}}=\eta^{\natural\ang3}$ we see
\[\sqrt[3]{\eta^{\natural\bot}}\in\iM(\Gamma_0(6)\cap\Gamma^0(2))_{\rho_3,\frac13,J[\natural]},\]
\[\sqrt[3]{\eta^{\natural\top}}\in\iM(\Gamma_0(2)\cap\Gamma^0(6))_{\rho_3,\frac13,J[\natural]}.\]

In particular $\sqrt[3]{\eta^{\natural\bot2}\eta^{\natural\top}}\in\cM(6)_1$, indeed
\[\sqrt[3]{\eta^{\natural\bot2}\eta^{\natural\top}}=\eta^{\natural\top}+4\sqrt{\eta^{\Diamond\top}}^{\ang{\frac13}}-4\sqrt{\eta^{\Diamond\top}}\]
since
\[\sqrt[3]{\eta^{\natural\bot2}\eta^{\natural\top}}\in1+4q^{\frac16}+4q^{\frac13}-2q^{\frac12}-4q^{\frac23}+4q+\C[[q]]q^{\frac76},\]
\[\eta^{\natural\top}\in1-2q^{\frac12}+4q+\C[[q]]q^{\frac32},\]
\begin{align*}\sqrt{\eta^{\Diamond\top}}&=\big(\sqrt{\eta^{\natural\top}}+2\sqrt{\eta^{\sharp\top}}\big)\sqrt{\eta^{\sharp\top}}=\eta^{\sharp\top\ang{\frac12}}+2\eta^{\sharp\top}\\&\in q^{\frac12}+q+q^{\frac32}-q^2+q^3+\C[[q]]q^{\frac72}.\end{align*}
Since $\sqrt[3]{\eta^{\natural\bot2}\eta^{\natural\top}}|_1(\BM 0&-1\\1&0\EM)=\frac{\sqrt3}{i}\sqrt[3]{\eta^{\natural\bot}\eta^{\natural\top2}}$ and $\eta^{\natural\top}|_1(\BM 0&-1\\1&0\EM)=\frac1{\sqrt3^3i}\eta^{\natural\bot}$, we have
\[\sqrt[3]{\eta^{\natural\bot}\eta^{\natural\top2}}=\z\frac19\big(\eta^{\natural\bot}+12\sqrt{\eta^{\Diamond\bot\ang3}}-4\sqrt{\eta^{\Diamond\bot}}\big).\]

\begin{Lem}\label{lev6}$\sqrt[3]{\eta^{\natural\bot}}=\sqrt[3]{\eta^{\natural\top}}+2\sqrt[6]{\eta^{\Diamond\top}}$ and $\sqrt[3]{\eta^{\natural\nwarrow}}=\sqrt[3]{\eta^{\natural\top}}+2\overline\omega\sqrt[6]{\eta^{\Diamond\top}}$.
\end{Lem}
\begin{proof}The former identity follows from
\begin{align*}&\z\frac14\big(\sqrt[3]{\eta^{\natural\bot}}-\sqrt[3]{\eta^{\natural\top}}\big)^3-2\sqrt{\eta^{\Diamond\top}}\\
&=\z\frac13\eta^{\natural\bot}+\sqrt{\eta^{\Diamond\bot\ang3}}-\frac13\sqrt{\eta^{\Diamond\bot}}-3\sqrt{\eta^{\Diamond\top}}^{\ang{\frac13}}+\sqrt{\eta^{\Diamond\top}}-\eta^{\natural\top}\\
&=\z\frac13\eta^{\natural\bot}-\frac13\sqrt{\eta^{\Diamond\bot}}-3\sqrt{\eta^{\Diamond\top}}^{\ang{\frac13}}\\
&=0\end{align*}
Acting $(\BM 1&1\\0&1\EM)$, the latter one follows from $\sqrt[3]{\eta^{\natural\top}}\in\C[[q^{\frac12}]]$ and $\sqrt[6]{\eta^{\Diamond\top}}\in \C[[q^{\frac12}]]q^{\frac16}$.

\end{proof}

Note $\sqrt[3]{\eta^{\natural8}}=\sqrt[3]{\eta^{\natural\bot}\eta^{\natural\top}\eta^{\natural\nwarrow}\eta^{\natural\swarrow}}$.

\begin{Thm}$\natural:\C\Big[\sqrt[3]{\eta^{\natural\bot}},\sqrt[3]{\eta^{\natural\top}},\sqrt[3]{\eta^{\natural4}}\Big]^{[\frac13,\frac13,\frac23]}\Big|_\M\bto\cM(6)$ and $\cM(6)\eta^4=\cS(6)$.\end{Thm}
\begin{proof}The dimension formula says for $k\in\M$
\[\dim\cM(6)_k=6k,\6\dim\cS(6)_k=\cM(6)_{k-2}.\]
We calculate
\[{\rm Dim\,LHS}=\frac{1+t^{\frac23}}{(1-t^{\frac13})^2}=1+\VS{k\in\N}2kt^{\frac k3}.\]
\end{proof}

Note $\sqrt[3]{\eta^{\Diamond4}}=\sqrt[6]{\eta^{\Diamond\top}}\sqrt{\eta^{\Diamond\bot\ang3}}$ and $\eta^4=\sqrt[3]{\eta^{\natural4}\eta^{\Diamond4}}$.

Also note
\begin{align*}\iE_{\rho_3}^{\ang{\frac12}}&=(\eta^{\flat\top\ang2}+2\eta^{\sharp\top\ang{\frac12}})+4(\eta^{\sharp\top\ang{\frac12}}+\eta^{\sharp\top})\\
&=\z\sqrt{\eta^{\natural\bot\ang3}\eta^{\natural\top}}+4\sqrt{\eta^{\sharp\bot\ang3}\eta^{\sharp\top}}\\
&=\eta^{\natural}\eta^{\natural\ang3}+4\eta^{\sharp}\eta^{\sharp\ang3}.\end{align*}
Since $x^2+3y^2=\big|(x+y)+\omega(x-y)\big|^2$, expanding out the above result gives for $n\in\N$
\[\#\big\{(a,b)\in\Z\times\Z \,\big|\, a^2+ab+b^2=n\big\}=6\VS{d|n}\rho_3(d).\]\\

\newpage

Formally we construct $S(6)^*=\big\{a+b \,\big|\, a\in\{0,\flat,\sharp\},\;b\in\{0,\top\}\big\}$. Let $S(6)^*$ be additive group $\simeq \Z/3\times\Z/2$ and
\[\fM(6)^*=\C\Big[\sqrt[6]{\eta^{\natural\top}},\sqrt[6]{\eta^{\flat\top}},\sqrt[6]{\eta^{\sharp\top}}\Big]^{[(\top,\frac16),(\flat+\top,\frac16),(\sharp+\top,\frac16)]}\]
be $\big(S(6)^*\times\frac16\M\big)$-graded.

With $e^{x}=t^{(x,0)}$ and $\tau^{\frac16}=t^{(0,\frac16)}$ we calculate
\[{\rm Dim}\,\fM(6)^*=\frac{(1-e^{\top}\tau^{\frac12})}{(1-e^{\top}\tau^{\frac16})(1-e^{\flat+\top}\tau^{\frac16})(1-e^{\sharp+\top}\tau^{\frac16})}\]
\[D_0={\rm Dim}\,\fM(6)^*|_{\{0,\flat,\sharp\}\times\frac13\M}=\dfrac{1+(2+e^{\flat}+e^{\sharp})\tau^{\frac13}+\tau^{\frac23}}{(1-e^{\flat}\tau^{\frac13})(1-e^{\sharp}\tau^{\frac13})},\]
\begin{align*}{\rm Dim}\,\fM(6)^*|_{\{0\}\times\frac13\M}&=\frac{1+2\tau^{\frac13}+4\tau^{\frac23}+4\tau+4\tau^{\frac43}+2\tau^{\frac53}+\tau^2}{(1-\tau)^2}\\
&=\frac{1+\tau^{\frac23}}{(1-\tau^{\frac13})^2}\end{align*}
In particular
\[\natural:\C\Big[\sqrt[3]{\eta^{\natural\bot}},\sqrt[3]{\eta^{\natural\top}},\sqrt[3]{\eta^{\natural4}}\Big]^{[\frac13,\frac13,\frac23]}\bto\fM(6)^*|_{\{0\}\times\frac13\M}.\]\\

By the former identity of Lem \ref{lev6}, naturally
\[\fM(6)^*\oplus\fM(6)^*\sqrt[6]{\eta^{\natural\bot}}\]
is a ring. Now, let
\[\fM(6)=\fM(6)^*\oplus\fM(6)^*\sqrt[6]{\eta^{\natural\bot}}\;\oplus\]
\[\6\C\Big[\sqrt[6]{\eta^{\natural\top}}\Big]\sqrt[6]{\eta^{\flat\bot}\eta^{\sharp\top}}\oplus\C\Big[\sqrt[6]{\eta^{\natural\top}}\Big]\sqrt[6]{\eta^{\flat\bot}}\oplus\C\Big[\sqrt[6]{\eta^{\natural\top}}\Big]\sqrt[6]{\eta^{\sharp\bot}}.\]
\begin{Lem}Naturally $\fM(6)$ is a ring.\end{Lem}

\begin{proof}
Acting $(\BM 1&3\\0&1\EM)$ on the former identity of Lem \ref{lev6}, we get
\[\sqrt[3]{\eta^{\flat\bot}}=\sqrt[3]{\eta^{\flat\top}}-2\sqrt[6]{\eta^{\natural\top}\eta^{\sharp\top}}\]
and
\begin{align*}\sqrt[6]{\eta^{\natural\bot}\eta^{\flat\bot}}&=\sqrt{\sqrt[3]{\eta^{\natural\top}\eta^{\flat\top}}-4\sqrt[6]{\eta^{\natural\top}\eta^{\flat\top}}\sqrt[3]{\eta^{\sharp\top}}+2\sqrt[6]{\eta^{\sharp\top}}\big(\sqrt{\eta^{\flat\top}}-\sqrt{\eta^{\natural\top}}\big)}\\
&=\sqrt[6]{\eta^{\natural\top}\eta^{\flat\top}}-2\sqrt[3]{\eta^{\sharp\top}}.\end{align*}
Hence $\fM(6)^*\oplus\fM(6)^*\sqrt[6]{\eta^{\natural\bot}}\oplus\fM(6)^*\sqrt[6]{\eta^{\flat\bot}}$ is a ring.

We see
\[\sqrt[3]{\eta^{\sharp\bot}}=\sqrt[3]{\eta^{\sharp\top}}+\sqrt[6]{\eta^{\natural\top}\eta^{\flat\top}}\]
thus
\begin{align*}\sqrt[6]{\eta^{\Diamond\bot}}&=\sqrt{\sqrt[6]{\eta^{\natural\top}}\big(\sqrt{\eta^{\flat\top}}-2\sqrt{\eta^{\sharp\top}}\big)-2\sqrt[3]{\eta^{\natural\top}}\sqrt[6]{\eta^{\Diamond\top}}+\sqrt[3]{\eta^{\Diamond\top}}}\\
&=\sqrt[3]{\eta^{\natural\top}}-\sqrt[6]{\eta^{\Diamond\top}},\end{align*}
acting $(\BM 1&3\\0&1\EM)$ we get
\[\sqrt[6]{\eta^{\natural\bot}\eta^{\sharp\bot}}=\sqrt[3]{\eta^{\flat\top}}+\sqrt[6]{\eta^{\natural\top}\eta^{\sharp\top}}.\]
Hence $\fM(6)^*\oplus\fM(6)^*\sqrt[6]{\eta^{\natural\bot}}\oplus\fM(6)^*\sqrt[6]{\eta^{\flat\bot}}\oplus\fM(6)^*\sqrt[6]{\eta^{\sharp\bot}}$ is a ring.

Remark
\[\sqrt[6]{\eta^{\natural\bot}\eta^{\natural\top}}=\sqrt[6]{\eta^{\sharp\bot}\eta^{\flat\top}}+\sqrt[6]{\eta^{\flat\bot}\eta^{\sharp\top}}\]
\[\sqrt[6]{\eta^{\flat\bot}\eta^{\flat\top}}=\sqrt[6]{\eta^{\sharp\bot}\eta^{\natural\top}}-\sqrt[6]{\eta^{\natural\bot}\eta^{\sharp\top}}\]
\[\sqrt[6]{\eta^{\sharp\bot}\eta^{\sharp\top}}=\z\frac12\big(\sqrt[6]{\eta^{\flat\bot}\eta^{\natural\top}}-\sqrt[6]{\eta^{\natural\bot}\eta^{\flat\top}}\big)\]
Hence 
\[\fM(6)=\fM(6)^*\oplus\fM(6)^*\sqrt[6]{\eta^{\natural\bot}}\oplus\fM(6)^*\sqrt[6]{\eta^{\flat\bot}\eta^{\sharp\top}}\oplus\C\Big[\sqrt[6]{\eta^{\natural\top}}\Big]\sqrt[6]{\eta^{\flat\bot}}\oplus\C\Big[\sqrt[6]{\eta^{\natural\top}}\Big]\sqrt[6]{\eta^{\sharp\bot}}\]
is also a ring.

At last, remark
\[\sqrt[6]{\eta^{\flat\top}}\sqrt[6]{\eta^{\flat\bot}\eta^{\sharp\top}}=\z\sqrt[6]{\eta^{\flat\bot}}\frac12\big(\sqrt[3]{\eta^{\natural\bot}}-\sqrt[3]{\eta^{\natural\top}}\big)\]
\[\sqrt[6]{\eta^{\sharp\top}}\sqrt[6]{\eta^{\flat\bot}\eta^{\sharp\top}}=\sqrt[6]{\eta^{\sharp\top}}\big(\sqrt[6]{\eta^{\natural\bot}\eta^{\natural\top}}-\sqrt[6]{\eta^{\sharp\bot}\eta^{\flat\top}}\big)\]
\end{proof}

Formally we construct $S(6)=\big\{a+b \,\big|\, a\in S(6)^*,\;b\in\{0,\bot\}\big\}$. Let $S(6)$ be additive group $\simeq S(6)^*\times\Z/2$. Moreover, let $\fM(6)$ be $\big(S(6)\times\frac16\M\big)$-graded by
\[\Big[\sqrt[6]{\eta^{\natural\bot}},\sqrt[6]{\eta^{\flat\bot}},\sqrt[6]{\eta^{\sharp\bot}}\Big]^{[(\bot,\frac16),(\flat+\bot,\frac16),(\sharp+\bot,\frac16)]}.\]

Note
\[{\rm Dim}\,\fM(6)={\rm Dim}\,\fM(6)^*\cdot(1+e^{\bot}\tau^{\frac16})+\frac{e^{\bot+\top}\tau^{\frac13}+e^{\flat+\bot}\tau^{\frac16}+e^{\sharp+\bot}\tau^{\frac16}}{1-e^{\top}\tau^{\frac16}}\]
and ${\rm Dim}\,\fM(6)|_{\{0,\flat,\sharp\}\times\frac13\M}=D_0$.

\newpage

\subsection{The case $N=12,24,48$}

\begin{Thm}Let
\[R=\fM(6)\Big[\sqrt[3]{\eta^{\natural2}},\sqrt[3]{\eta^{\flat2}},\sqrt[3]{\eta^{\sharp2}}\Big]^{[(0,\frac13),(\sharp,\frac13),(\flat,\frac13)]}\]
then $\natural:R|_{\{0\}\times\M}\bto\cM(12)$ and $\cM(12)\eta^2=\cS(12)$.\end{Thm}
\begin{proof}The dimension formula says for $k\in\N$
\[\dim\cM(12)_k=48k-24+\delta_{k,1},\6\dim\cS(12)_k=\dim\cM(12)_{k-1}.\]

With $e^x=t^{(x,0)}$ and $\tau^{\frac16}=t^{(0,\frac16)}$ we calculate
\[{\rm Dim}\,R|_{\{0,\flat,\sharp\}\times\frac13\M}=D_1\]
where $D_1=D_0(1+\tau^{\frac13})(1+e^{\sharp}\tau^{\frac13})(1+e^{\flat}\tau^{\frac13})$,
\[{\rm Dim}\,R|_{\{0\}\times\frac13\M}=\frac{(1+\tau^{\frac13})(1+6\tau^{\frac23}+\tau^{\frac43})}{(1-\tau^{\frac13})^2}.\]
\end{proof}

\begin{Lem}Let
\[R=\fM(6)\Big[\sqrt[3]{\eta^{\natural}},\sqrt[3]{\eta^{\flat}},\sqrt[3]{\eta^{\sharp}}\Big]^{[(0,\frac16),(\flat,\frac16),(\sharp,\frac16)]}\]
then $\natural:R|_{\{0\}\times\frac12\M}\bto\iM(24)_{\ang{\xi_8,\rho_3},\frac12\M}$.\end{Lem}
\begin{proof}The dimension formula says
\[\dim\iM(24)_{\ang{\rho_4,\rho_3},k}=192k-144+10\delta_{k,1} \text{ for }k\in\N\]
\begin{align*}\dim\iM(24)_{\ang{\rho_4,\rho_3},\frac k2}&\leqq\dim\iS(24)_{\ang{\rho_4,\rho_3},\frac{k+1}2}\\
&=(96k-144)^++10\delta_{k,1}+\delta_{k,3} \text{ for }k\in2\M+1\end{align*}

With $e^x=t^{(x,0)}$ and $\tau^{\frac16}=t^{(0,\frac16)}$ we calculate
\[{\rm Dim}\,R|_{\{0,\flat,\sharp\}\times\frac13\M}=D_2\]
where $D_2=D_1(1+\tau^{\frac16})(1+e^{\flat}\tau^{\frac16})(1+e^{\sharp}\tau^{\frac16})$
\[{\rm Dim}\,R|_{\{0\}\times\frac13\M}=\frac{(1+\tau^{\frac16}+\tau^{\frac13}+\tau^{\frac12})(1-2\tau^{\frac16}+4\tau^{\frac13}+2\tau^{\frac23}+4\tau-2\tau^{\frac76}+\tau^{\frac43})}{(1-\tau^{\frac16})^2}.\]
\end{proof}

We easily see
\[\natural:\C\Big[\sqrt[3]{\eta^{\natural\bot}},\sqrt[3]{\eta^{\natural\top}},\sqrt[3]{\eta^{\natural}}\Big]^{[\frac13,\frac13,\frac16]}\Big|_{\frac12\M}\bto\iM(6)_{\ang{\xi_8,\rho_3}}.\]

Note
\[\sqrt[3]{\eta^{\natural}}\sqrt[6]{\eta^{\Diamond\top}}=\z\frac12(\eta^{\natural\ang{\frac13}}-\eta^{\natural\ang3})=\theta_{{\tt 1}_3}^{\ang{\frac16}},\]
\[\sqrt[3]{\eta^{\natural}}\sqrt[6]{\eta^{\Diamond\bot}}=\z\frac12(3\eta^{\natural\ang3}-\eta^{\natural\ang{\frac13}})=\big(\VS{n\in\Z}\omega^nq^{n^2}\big)^{\ang{\frac16}}.\]
Acting $(\BM 1&3\\0&1\EM)$ on these identities, we get $\sqrt[3]{\eta^{\flat}}\sqrt[6]{\eta^{\natural\top}\eta^{\sharp\top}}=\big(\VS{n\in\Z}(-1)^n{\tt 1}_3(n)q^{n^2}\big)^{\ang{\frac16}}$ and $\sqrt[3]{\eta^{\flat}}\sqrt[6]{\eta^{\natural\bot}\eta^{\sharp\bot}}=\big(\VS{n\in\Z}(-\omega)^nq^{n^2}\big)^{\ang{\frac16}}$.

We also see
\[\sqrt[3]{\eta^{\sharp}}\sqrt[6]{\eta^{\natural\top}\eta^{\flat\top}}=\eta^{\sharp\ang{\frac13}}-\eta^{\sharp\ang3}=\theta_{{\tt 1}_6}^{\ang{\frac1{24}}},\]
\[\sqrt[3]{\eta^{\sharp}}\sqrt[6]{\eta^{\natural\bot}\eta^{\flat\bot}}=\eta^{\sharp\ang{\frac13}}-3\eta^{\sharp\ang3}=\big(\VS{n\in\Z}(-\omega)^n{\tt 1}_2(n)q^{n^2}\big)^{\ang{\frac1{24}}}.\]\\

Moreover we have $\dim\iS(6)_{\xi_8^3,\frac32}=1$ and $\eta\eta^{\Diamond}=\sqrt[3]{\eta^{\natural}\eta^{\Diamond4}}=\Theta_{\rho_3}^{\ang{\frac16}}$. Acting $(\BM 1&3\\0&1\EM)$, we get $\eta\eta^{\natural}\eta^{\sharp}=\big(\VS{n\in\N}(-1)^{n-1}\rho_3(n)nq^{n^2}\big)^{\ang{\frac16}}$ which is a classical result due to K\"ohler-Macdonald. Lemma \ref{rel21}, \ref{rel22} lead to
\[(\eta\eta^{\Diamond})^{\ang{\frac14}}=\sqrt[3]{(\eta^{\natural}+2\eta^{\sharp})(\eta^{\natural2}+4\eta^{\sharp2})\eta^{\natural}\eta^{\sharp}(\eta^{\natural}-2\eta^{\sharp})^4}=\eta\eta^{\flat}(\eta^{\natural}-2\eta^{\sharp})\]
thus
\[\eta\eta^{\natural}\eta^{\flat}=(\eta\eta^{\Diamond})^{\ang{\frac14}}+2\eta\eta^{\Diamond}=\Theta_{{\tt 1}_2\rho_3}^{\ang{\frac1{24}}}.\]\\

In addition, let
\[R=\fM(6)\Big[\sqrt[6]{\eta^{\natural}},\sqrt[6]{\eta^{\flat}},\sqrt[6]{\eta^{\sharp}}\Big]^{[(\frac1{12}),(\sharp,\frac1{12}),(\flat,\frac1{12})]}\]
then $\natural:R|_{\{0\}\times\frac12\M}\to\iM(48)_{\ang{\xi_8,\rho_8,\rho_3},\frac12\M}$.

Note $\sqrt[3]{\eta^{\flat\top}}\sqrt[6]{\eta^{\natural}\eta^{\sharp}}=\big(\sqrt[3]{\eta^{\natural}}\sqrt[6]{\eta^{\Diamond\top}}\big)^{\ang{\frac12}}=\theta_{{\tt 1}_6}^{\ang{\frac1{48}}}$. Acting $(\BM 1&1\\0&1\EM)$, $\sqrt[3]{\eta^{\natural\top}}\sqrt[6]{\eta^{\Diamond}}=\theta_{{\tt 1}_3\rho_8}^{\ang{\frac1{48}}}$. In addition
\[\sqrt[3]{\eta^{\natural\bot}}\sqrt[6]{\eta^{\Diamond}}=(\theta_{{\tt 1}_3\rho_8}+2\theta_{\rho_8}^{\ang9})^{\ang{\frac1{48}}}=(\theta_{\rho_8}+3\theta_{\rho_8}^{\ang9})^{\ang{\frac1{48}}}.\]

Also $\sqrt[3]{\eta^{\flat2}}\sqrt[6]{\eta^{\natural}\eta^{\sharp}}=\eta^{\ang{\frac12}}=\theta_{\rho_4\rho_3}^{\ang{\frac1{48}}}$. Acting $(\BM 1&1\\0&1\EM)$, $\sqrt{\eta\eta^{\natural}}=\sqrt[3]{\eta^{\natural2}}\sqrt[6]{\eta^{\Diamond}}=\theta_{\rho_4\rho_8\rho_3}^{\ang{\frac1{48}}}$.

At last, we have $\sqrt{\eta\eta^{\natural5}}=\Theta_{\rho_8\rho_3}^{\ang{\frac1{48}}}$.

\newpage

\subsection{The case $N=18$}

Let
\[\fM(6)'=\fM(6)\oplus\fM(6)\sqrt[6]{\eta^{\natural\nwarrow}}\]
\[\6\oplus\C\Big[\sqrt[6]{\eta^{\natural\top}}\Big]\sqrt[6]{\eta^{\flat\top}\eta^{\sharp\nwarrow}}\oplus\C\Big[\sqrt[6]{\eta^{\natural\top}}\Big]\sqrt[6]{\eta^{\natural\swarrow}}\]
\[\6\oplus\C\Big[\sqrt[6]{\eta^{\natural\top}}\Big]\sqrt[6]{\eta^{\flat\nwarrow}}\oplus\C\Big[\sqrt[6]{\eta^{\natural\top}}\Big]\sqrt[6]{\eta^{\flat\swarrow}}\oplus\C\Big[\sqrt[6]{\eta^{\natural\top}}\Big]\sqrt[6]{\eta^{\sharp\nwarrow}}\oplus\C\Big[\sqrt[6]{\eta^{\natural\top}}\Big]\sqrt[6]{\eta^{\sharp\swarrow}}.\]

Formally we construct $S(6)'=\big\{a+b \,\big|\, a\in S(6),\;b\in\{0,\nwarrow\}\big\}$. Let $S(6)'$ be additive group $\simeq S(6)\times\Z/2$ and ${\swarrow}=\bot+\top+{\nwarrow}$. Moreover, let $\fM(6)'$ be $\big(S(6)'\times\frac16\M\big)$-graded by
\[\Big[\sqrt[6]{\eta^{\natural\nwarrow}},\sqrt[6]{\eta^{\flat\nwarrow}},\sqrt[6]{\eta^{\sharp\nwarrow}},\sqrt[6]{\eta^{\natural\swarrow}},\sqrt[6]{\eta^{\flat\swarrow}},\sqrt[6]{\eta^{\sharp\swarrow}}\Big]\]
\[\6^{[(\nwarrow,\frac16),(\flat+\nwarrow,\frac16),(\sharp+\nwarrow,\frac16),(\swarrow,\frac16),(\flat+\swarrow,\frac16),(\sharp+\swarrow,\frac16)]}.\]

Note $\sqrt{\eta^{\top}}=\sqrt[6]{\eta^{\natural\top}\eta^{\flat\top}\eta^{\sharp\top}}$ has weight $(\top,\frac12)$ in $\fM(6)$. Since
\[\sqrt[6]{\eta^{\top}\eta^{\bot2}}\sqrt{\eta^{\top}}=\sqrt[3]{\eta^{\top2}\eta^{\bot}}\in\iM(9)_{\rho_3,1},\]
\[\sqrt[6]{\eta^{\bot}\eta^{\bot\top}}\sqrt{\eta^{\top}}=\eta\eta^{\ang3}\in\iM(\Gamma_0(18)\cap\Gamma^0(6))_{\rho_3,1}\]
$\sqrt[6]{\eta^{\top}\eta^{\bot2}},\sqrt[6]{\eta^{\bot}\eta^{\bot\top}}$ also should have weight $(\top,\frac12)$.

\begin{Thm}Let
\[R=\fM(6)'|_{S(6)'\times\frac12\M}\Big[\sqrt[6]{\eta^{\top}},\sqrt[6]{\eta^{\bot}},\sqrt[6]{\eta^{\nwarrow}},\sqrt[6]{\eta^{\swarrow}}\Big]^{[(\top,\frac16),(\bot,\frac16),(\nwarrow,\frac16),(\swarrow,\frac16)]}\]
then $\natural:R|_{\{0\}\times\frac13\M}\bto\cM(18)_{\frac13\M}$. And $\cM(18)_{\frac13\M}\eta^{\frac43}=\cS(18)_{\frac13\M}$.\end{Thm}
\begin{proof}The dimension formula says
\[\dim\cM(18)_k=162k-108+4\delta_{k,1} \text{ for }k\in\N\]
We see
\[R_{S(6)'\times\{\kappa\}}=\fM(6)'|_{S(6)'\times\frac12\M}\oplus\VT{f\in A}\fM(6)'|_{S(6)'\times\{\kappa-\frac12\}}f\oplus\VT{f\in B}\fM(6)'|_{S(6)'\times\{\kappa-1\}}f\]
where
\[A=\big\{\sqrt[6]{\eta^{\top2}\eta^{\bot}},\sqrt[6]{\eta^{\top}\eta^{\bot2}},\sqrt[6]{\eta^{\top}\eta^{\nwarrow}\eta^{\swarrow}},\sqrt[6]{\eta^{\bot}\eta^{\nwarrow}\eta^{\swarrow}}\]
\[\6\sqrt[6]{\eta^{\top}\eta^{\nwarrow2}},\sqrt[6]{\eta^{\bot}\eta^{\nwarrow2}},\sqrt[6]{\eta^{\top}\eta^{\swarrow2}},\sqrt[6]{\eta^{\bot}\eta^{\swarrow2}},\sqrt[6]{\eta^{\top}\eta^{\bot}\eta^{\nwarrow}},\sqrt[6]{\eta^{\top}\eta^{\bot}\eta^{\swarrow}},\]
\[\6\sqrt[6]{\eta^{\top2}\eta^{\nwarrow}},\sqrt[6]{\eta^{\bot2}\eta^{\nwarrow}},\sqrt[6]{\eta^{\top2}\eta^{\swarrow}},\sqrt[6]{\eta^{\bot2}\eta^{\swarrow}},\sqrt[6]{\eta^{\nwarrow2}\eta^{\swarrow}},\sqrt[6]{\eta^{\nwarrow}\eta^{\swarrow2}}\big\}\]
\[B=\big\{\sqrt[6]{\eta^{\top2}\eta^{\nwarrow2}\eta^{\swarrow2}},\sqrt[6]{\eta^{\bot2}\eta^{\nwarrow2}\eta^{\swarrow2}},\sqrt[6]{\eta^{\top2}\eta^{\bot2}\eta^{\nwarrow}\eta^{\swarrow}},\sqrt[6]{\eta^{\top}\eta^{\bot}\eta^{\nwarrow2}\eta^{\swarrow2}}\]
\[\6\sqrt[6]{\eta^{\top2}\eta^{\bot2}\eta^{\nwarrow2}},\sqrt[6]{\eta^{\top2}\eta^{\bot2}\eta^{\swarrow2}}\]
\[\6\sqrt[6]{\eta^{\top}\eta^{\nwarrow2}\eta^{\swarrow}},\sqrt[6]{\eta^{\bot}\eta^{\nwarrow2}\eta^{\swarrow}},\sqrt[6]{\eta^{\top}\eta^{\nwarrow}\eta^{\swarrow2}},\sqrt[6]{\eta^{\bot}\eta^{\nwarrow}\eta^{\swarrow2}}\big\}\]
thus $\dim\,R_{0,\frac k6}=k+16(k-3)+10(k-6)$.

The dimension formula also says
\[\z\dim\cS(18)_{k}=162(k-\frac23)-108+\delta_{k,2} \text{ for }k\geq2.\]
\end{proof}

\newpage

\subsection{The case $N=5$}

Define
\[u_5=\z\sqrt[5]{\frac12\big(\iE_{\chi_5}+\iE_{\overline{\chi_5}}\big)}.\]
\[a_5=\z\sqrt[5]{\frac i2\big(\iE_{\chi_5}-\iE_{\overline{\chi_5}}\big)}.\]

\begin{Lem}\label{lev3}$(u_5-ia_5)^5,(u_5-ia_5)(u_5+ia_5)^4\in\iM(5)_{\chi_5,1}$.\end{Lem}
\begin{proof}
Note $\iG_{\chi_5}=\frac1{3-i}\big(\iE_{\chi_5}-\iE_{\chi_5}^{\ang5}\big)$. We have $\dim\big\{f+\overline f\,\big|\,f\in\iM(5)_{\chi_5,5}\big\}=13$ and
\[\z u_5^4a_5=\frac14(\iG_{\chi_5}+\iG_{\rho_5|\overline{\chi_5}}+\iG_{\chi_5}^c + \iG_{\rho_5|\overline{\chi_5}}^c)^{\ang{\frac15}}.\]
Hence $\dfrac{a_5}{u_5}=\dfrac{a_5^5}{u_5^4a_5}$ is $\Gamma(5)$-invariant and so are
\[u_5^3a_5^2=u_5^4a_5\frac{a_5}{u_5},\6 u_5^2a_5^3=u_5^3a_5^2\frac{a_5}{u_5},\6 u_5a_5^4=u_5^2a_5^3\frac{a_5}{u_5}.\]

Moreover, $\dim\iM(5)_1=6$ and
\[(u_5-ia_5)^5=\z\frac{-3-i}4\iE_{\chi_5}^{\ang{\frac15}}+\frac{7+i}4\iE_{\chi_5}+(\frac52-5i)\iG_{\rho_5|\overline{\chi_5}}^{\ang{\frac15}}.\]
\[(u_5-ia_5)(u_5+ia_5)^4=\z\frac{-7+11i}{20}\iE_{\chi_5}^{\ang{\frac15}}+\frac{27-11i}{20}\iE_{\chi_5}+(\frac12+i)\iG_{\rho_5|\overline{\chi_5}}^{\ang{\frac15}},\]
\end{proof}
Now, let $J[5]=J\big[\chi_5,u_5-ia_5\big]$ on $\Gamma_0(5)\cap\Gamma^0(5)$ and $|_{\frac k5}=|_{\frac k5,J[5]}$ unless otherwise noted, then $u_5-ia_5\in\iM(5)_{\chi_5,\frac15}$ and $u_5+ia_5\in\iM(5)_{\overline{\chi_5},\frac15}$.

The dimension formula says $\dim\iM(\Gamma_0(5))_{\rho_5,2}=2$ and $u_5a_5=\sqrt[5]{\eta^{"5"}}$.
Also $\dim\iM(5)_{\rho_5,2}=6$ and $u_5^2-a_5^2-u_5a_5=\sqrt[5]{\eta^{/5/}}$.

Note $\eta^{\frac{24}5}=\eta^{/5/\ang5}\sqrt[5]{\eta^{"5"}}=\eta^{"5"\ang{\frac15}}\sqrt[5]{\eta^{/5/}}\in\cM(5)_{\frac{12}5}$.\\

\begin{Thm}$\natural:\C\big[u_5,a_5\big]^{[\frac15,\frac15]}\bto\cM(5)_{\frac15\M}$ and $\cM(5)_{\frac15\M}\eta^{\frac{24}5}=\cS(5)_{\frac15\M}$.\end{Thm}
\begin{proof}The dimension formula says for $k\in\M$
\[\z\dim\cM(5)_k=5k+1,\6\dim\cS(5)_k=\big(5(k-\frac{12}5)+1\big)^+.\]

We easily $\dim\cM(5)_{\frac{k-1}5}\leqq\dim\cM(5)_{\frac k5}-1$ thus
\[{\rm Dim}\,\cM(5)_{\frac15\M}=\VS{k\in\M}(k+1)t^{\frac k5}.\]
\end{proof}

Indeed, the generators are essentially the famous Rogers-Ramanujan functions;
\[u_5=\VP{n\in\N}(1-q^n)^{\frac25}\VP{n\in\M}\dfrac1{(1-q^{5n+1})(1-q^{5n+4})},\]
\[a_5=q^{\frac15}\VP{n\in\N}(1-q^n)^{\frac25}\VP{n\in\M}\dfrac1{(1-q^{5n+2})(1-q^{5n+3})}.\]

\newpage

\subsection{The case $N=10$}

First, put
\[\iE_{10}=i\iE_{\chi_5}+(1-i)\iE_{\chi_5}^{\ang2}.\]
This may be a modular form of weight $\frac13$ since the dimension formula says\\$\dim\cM(\Gamma_0(10))_{\chi_5,3}=5$ and $\sqrt[3]{\iE_{10}\iE_{10}^{c2}}=\frac{1-i}{3-i}\iE_{\chi_5}+\frac2{3-i}\iE_{\chi_5}^{\ang2}$. Define
\[\z u_{10}=\sqrt[5]{\frac12\big(\sqrt[3]{\iE_{10}}+\sqrt[3]{\iE_{10}^c}\big)},\]
\[\z a_{10}=\sqrt[5]{\frac1{2i}\big(\sqrt[3]{\iE_{10}}-\sqrt[3]{\iE_{10}^c}\big)},\]
\[\z u_{10}'=\sqrt[5]{u_{10}^5+a_{10}^5},\6\z u_{10}''=\sqrt[5]{u_{10}^5-a_{10}^5}\]
then
\[\z u_{10}^2u_{10}''=\sqrt[5]{\frac18\big(\sqrt[3]{\iE_{10}}+\sqrt[3]{\iE_{10}^c}\big)^2\big((1+i)\sqrt[3]{\iE_{10}}+(1-i)\sqrt[3]{\iE_{10}^c}\big)}=u_5^{\ang2}\]
\[\z a_{10}^2u_{10}'=\sqrt[5]{-\frac18\big(\sqrt[3]{\iE_{10}}-\sqrt[3]{\iE_{10}^c}\big)^2\big((1-i)\sqrt[3]{\iE_{10}}+(1+i)\sqrt[3]{\iE_{10}^c}\big)}=a_5^{\ang2}\]
and similarly $u_{10}u_{10}'^2=u_5$, $a_{10}u_{10}''^2=a_5$. In particular, $u_{10}a_{10}u_{10}'u_{10}''=\sqrt[3]{u_5a_5u_5^{\ang2}u_5^{\ang2}}$ and
\[u_{10}a_{10}=\frac{u_{10}^2u_{10}''a_{10}^2u_{10}'}{u_{10}a_{10}u_{10}'u_{10}''}=\sqrt[3]{u_5a_5}^{\sharp}\]
\[u_{10}'u_{10}''=\frac{u_{10}u_{10}'^2a_{10}u_{10}''^2}{u_{10}a_{10}u_{10}'u_{10}''}=\sqrt[3]{u_5a_5}^{\flat\ang2}\]

The dimension formula says $\dim\cM(\Gamma_1(10))_2=7$ and
\[\z u_{10}'^5u_{10}''^5-4u_{10}^5a_{10}^5=\sqrt[3]{\eta^{\flat/5/\ang{10}}},\]
\[\z u_{10}'^5u_{10}''^5+u_{10}^5a_{10}^5=\sqrt[3]{\eta^{\sharp/5/\ang5}}.\]

Note $\eta^{\sharp\frac45}=\sqrt{u_{10}a_{10}(u_{10}'^5u_{10}''^5+u_{10}^5a_{10}^5)}$. We have
\[u_{10}a_{10}^3u_{10}'^3u_{10}''^2\eta^{\sharp\frac45}=(u_5^3a_5^2)^{\ang{\frac12}}-u_5^4a_5+(u_5^2a_5^3)^{\ang2}\]
\[u_{10}^3a_{10}u_{10}'^2u_{10}''^3\eta^{\sharp\frac45}=(u_5^2a_5^3)^{\ang{\frac12}}+u_5a_5^4-(u_5^3a_5^2)^{\ang2}\]
and
\[u_{10}^4u_{10}'^4u_{10}''\eta^{\sharp\frac45}=\frac{(u_{10}a_{10}^3u_{10}'^3u_{10}''^2\eta^{\sharp\frac45})^2}{u_{10}^3a_{10}u_{10}'^2u_{10}''^3\eta^{\sharp\frac45}}\frac{u_{10}^5}{a_{10}^5}\]
\[a_{10}^4u_{10}'u_{10}''^4\eta^{\sharp\frac45}=\frac{(u_{10}a_{10}^3u_{10}'^3u_{10}''^2\eta^{\sharp\frac45})(u_{10}^3a_{10}u_{10}'^2u_{10}''^3\eta^{\sharp\frac45})}{u_{10}^4u_{10}'^4u_{10}''\eta^{\sharp\frac45}}\]
\[u_{10}^2a_{10}^2u_{10}''^5\eta^{\sharp\frac45}=\frac{(u_{10}^3a_{10}u_{10}'^2u_{10}''^3\eta^{\sharp\frac45})^2}{u_{10}^4u_{10}'^4u_{10}''\eta^{\sharp\frac45}}\]
\[u_{10}^2a_{10}^2u_{10}'^5\eta^{\sharp\frac45}=u_{10}^2a_{10}^2u_{10}''^5\eta^{\sharp\frac45}\frac{u_5'^5}{u_5''^5}\]\\

Now, let
\[\cM(10)^*=\C\big[u_{10},a_{10},u_{10}',u_{10}'',\eta^{\sharp\frac45}\big]^{[(\frac15,0,\frac1{15}),(0,\frac15,\frac1{15}),(\frac25,0,\frac1{15}),(0,\frac25,\frac1{15}),(\frac35,\frac35,\frac25)]}\]
be $(\frac15\Z/\Z\times\frac15\Z/\Z\times\frac1{15}\M)$-graded.

\newpage

\begin{Thm}$\natural:\cM(10)^*|_{\{0\}\times\{0\}\times\M}\bto\cM(10)$.
\end{Thm}
\begin{proof}With $e^{\frac15}=t^{(\frac15,0,0)}$, $f^{\frac15}=t^{(0,\frac15,0)}$ and $\tau^{\frac1{15}}=t^{(0,0,\frac1{15})}$ we calculate
\[{\rm Dim}\,\cM(10)^*=\frac{(1-\tau^{\frac13})^2(1+e^{\frac35}f^{\frac35}\tau^{\frac25})}{(1-e^{\frac15}\tau^{\frac1{15}})(1-e^{\frac25}\tau^{\frac1{15}})(1-f^{\frac15}\tau^{\frac1{15}})(1-f^{\frac25}\tau^{\frac1{15}})}\]
\begin{align*}&{\rm Dim}\,\cM(10)^*|_{\{0\}\times\frac15\Z/\Z\times\frac1{15}\M}=\frac{(1-\tau^{\frac25})^2}{5(1-f^{\frac15}\tau^{\frac1{15}})(1-f^{\frac25}\tau^{\frac1{15}})}\Big(\frac{1+e^{\frac35}f^{\frac35}\tau^{\frac25}}{(1-e^{\frac15}\tau^{\frac1{15}})(1-e^{\frac25}\tau^{\frac1{15}})}\\
&\6+\frac{1+\zeta_5^3e^{\frac35}f^{\frac35}\tau^{\frac25}}{(1-\zeta_5e^{\frac15}\tau^{\frac1{15}})(1-\zeta_5^2e^{\frac25}\tau^{\frac1{15}})}+\frac{1+\zeta_5e^{\frac35}f^{\frac35}\tau^{\frac25}}{(1-\zeta_5^2e^{\frac15}\tau^{\frac1{15}})(1-\overline{\zeta_5}e^{\frac25}\tau^{\frac1{15}})}\\
&\6+\frac{1+\overline{\zeta_5}e^{\frac35}f^{\frac35}\tau^{\frac25}}{(1-\zeta_5^3e^{\frac15}\tau^{\frac1{15}})(1-\zeta_5e^{\frac25}\tau^{\frac1{15}})}+\frac{1+\zeta_5^2e^{\frac35}f^{\frac35}\tau^{\frac25}}{(1-\overline{\zeta_5}e^{\frac15}\tau^{\frac1{15}})(1-\zeta_5^3e^{\frac25}\tau^{\frac1{15}})}\Big)\\
&=\frac{1+\tau^{\frac3{15}}+\tau^{\frac4{15}}+\tau^{\frac6{15}}+\tau^{\frac7{15}}+f^{\frac35}\tau^{\frac25}(\tau^{\frac1{15}}+\tau^{\frac2{15}}+\tau^{\frac4{15}}+\tau^{\frac5{15}}+\tau^{\frac8{15}})}{(1-f^{\frac15}\tau^{\frac1{15}})(1-f^{\frac25}\tau^{\frac1{15}})}\end{align*}
\begin{align*}&{\rm Dim}\,\cM(10)^*|_{\{0\}\times\{0\}\times\frac1{15}\M}\\
&=\frac{(1+\tau^{\frac3{15}}+\tau^{\frac4{15}}+\tau^{\frac6{15}}+\tau^{\frac7{15}})^2+\tau^{\frac25}(\tau^{\frac1{15}}+\tau^{\frac2{15}}+\tau^{\frac4{15}}+\tau^{\frac5{15}}+\tau^{\frac8{15}})^2}{(1-\tau^{\frac13})^2}\end{align*}
We calculate
\[\frac{1+5t^{\frac23}+4t}{(1-t^{\frac13})^2}=\VS{k\in\M}((k+1)+5(k-1)^++4(k-2)^+)t^{\frac k3}=1+4t^{\frac13}+\VS{k\in\N}(10k-12)t^{\frac k3}.\]

The dimension formula says $\dim\cM(10)_k=30k-12$ for $k\in\N$.
\end{proof}

\vspace{0.5cm}

The dimension formula also says $\dim\cS(10)_k=\delta_{k,2}+(30(k-\frac65)-12)^+$.

Note $\eta^{\flat\frac85}=u_{10}'u_{10}''(u_{10}'^5u_{10}''^5-4u_{10}^5a_{10}^5)$, $\eta^{\frac{12}5}=\eta^{\sharp\frac45}\eta^{\flat\frac85}$ and
\[\dim\cS(10)=\dfrac{4t^{\frac7{15}}+5t^{\frac45}+t^{\frac{22}{15}}}{(1-t^\frac13)^2}t^{\frac65}\Big|_\M\]
thus
\[\natural:\cM(10)^*\eta^{\frac{12}5}\big|_\M\bto\cS(10).\]

\newpage

\subsection{The case $N=7$}

Define
\[\z u_7=\sqrt[7]{\frac1{\overline\omega-\omega}\big(\overline\omega\sqrt{\iE_{\chi_7}}-\omega\sqrt{\iE_{\overline{\chi_7}}}\big)},\]
\[\z a_7=\sqrt[7]{\frac1{\overline\omega-\omega}\big(\sqrt{\iE_{\chi_7}}-\sqrt{\iE_{\overline{\chi_7}}}\big)},\]
\[u_7'=\sqrt[7]{u_7^7+a_7^7},\]
then
\[(u_7a_7u_7')^7=\z\frac1{(\omega-\overline\omega)^3}(\sqrt{\iE_{\chi_7}}^3-\sqrt{\iE_{\overline{\chi_7}}}^3)\]
and $(u_7a_7u_7')^{14}\in\iM(7)_{\rho_7,3}$. So, let $J[14]=J\big[\sqrt{\rho_7},\sqrt[3]{u_7a_7u_7'}\big]$ on $\Gamma_0(7)\cap\Gamma^0(7)$ and $|_{\frac k{14}}=|_{\frac k{14},J[14]}$ unless otherwise noted, then $u_7a_7u_7'\in\iM(7)_{\sqrt{\rho_7}^3,\frac3{14}}$. Since
\[(u_7a_7u_7')^7(u_7^7-\omega a_7^7)=\z\frac1{(\omega-\overline\omega)^3}(\iE_{\chi_7}^2-\iE_{\rho_7}\iE_{\overline{\chi_7}})\]
we get $u_7^7-\omega a_7^7\in\iM(7)_{\sqrt{\rho_7}\chi_7,\frac7{14}}$. Moreover
\[u_7'u_7^3+a_7u_7'^3-u_7a_7^3\in\iM(7)_{{\tt 1}_7,\frac4{14}},\]
\[u_7'u_7^3+\omega a_7u_7'^3-\overline\omega u_7a_7^3\in\iM(7)_{\chi_7,\frac4{14}}.\]
We see $u_7'^3u_7^2=\dfrac{(u_7'u_7^3)^3}{u_7^7}\in\cM(7)_{\frac5{14}}$ and similarly $u_7^3a_7^2,a_7^3u_7'^2\in\cM(7)_{\frac5{14}}$.\\

The dimension formula says $\dim\iM(\Gamma_0(7))_{{\tt 1}_7,6}=5$ and
\[\z \frac{11-13\zeta_6}9\sqrt{\iE_{\chi_7}}^3+\frac{11-13\overline{\zeta_6}}9\sqrt{\iE_{\overline{\chi_7}}}^3=\sqrt{\eta^{/7/\ang7}},\]
\[\z u_7a_7u_7'=\sqrt[14]{\eta^{"7"}}.\]
Also $\dim\iM(7)_{{\tt 1}_7,6}=28$ and
\[\z u_7u_7'^5-u_7^5a_7+u_7'a_7^5-u_7^2u_7'^2a_7^2=\sqrt[7]{\eta^{/7/}}.\]
Note $\eta^{\frac{24}7}=\sqrt{\eta^{/7/\ang7}}\sqrt[14]{\eta^{"7"}}\in\cM(7)^*_{\frac{12}7}$.

Now, let
\[\cM(7)^*=\C[u_7,a_7,u_7']^{[(\frac17,\frac1{14}),(\frac27,\frac1{14}),(\frac47,\frac1{14})]}\]
be $(\frac17\Z/\Z\times\frac1{14}\M)$-graded.

\begin{Thm}$\natural:\cM(7)^*|_{\{0\}\times\M}\bto\cM(7)$.
\end{Thm}
\begin{proof}With $e^{\frac17}=t^{(\frac17,0)}$ and $\tau^{\frac1{14}}=t^{(0,\frac1{14})}$ we calculate
\[{\rm Dim}\,\cM(7)^*=\frac{1-\tau^{\frac12}}{(1-e^{\frac17}\tau^{\frac1{14}})(1-e^{\frac27}\tau^{\frac1{14}})(1-e^{\frac47}\tau^{\frac1{14}})}\]
\begin{align*}&{\rm Dim}\,\cM(7)^*|_{\{0\}\times\frac1{14}\M}=\frac{1-\tau^{\frac12}}{1-\tau^{\frac1{14}}}\frac17\Big(\frac1{(1-e^{\frac17}\tau^{\frac1{14}})(1-e^{\frac27}\tau^{\frac1{14}})(1-e^{\frac47}\tau^{\frac1{14}})}\\
&\6+\frac1{(1-\zeta_7e^{\frac17}\tau^{\frac1{14}})(1-\zeta_7^2e^{\frac27}\tau^{\frac1{14}})(1-\zeta_7^4e^{\frac47}\tau^{\frac1{14}})}+\frac1{(1-\zeta_7^2e^{\frac17}\tau^{\frac1{14}})(1-\zeta_7^4e^{\frac27}\tau^{\frac1{14}})(1-\zeta_7e^{\frac47}\tau^{\frac1{14}})}\\
&\6+\frac1{(1-\zeta_7^3e^{\frac17}\tau^{\frac1{14}})(1-\overline{\zeta_7}e^{\frac27}\tau^{\frac1{14}})(1-\zeta_7^5e^{\frac57}\tau^{\frac1{14}})}+\frac1{(1-\zeta_7^4e^{\frac17}\tau^{\frac1{14}})(1-\zeta_7e^{\frac27}\tau^{\frac1{14}})(1-\zeta_7^2e^{\frac47}\tau^{\frac1{14}})}\\
&\6+\frac1{(1-\zeta_7^5e^{\frac17}\tau^{\frac1{14}})(1-\zeta_7^3e^{\frac27}\tau^{\frac1{14}})(1-\overline{\zeta_7}e^{\frac47}\tau^{\frac1{14}})}+\frac1{(1-\overline{\zeta_7}e^{\frac17}\tau^{\frac1{14}})(1-\zeta_7^5e^{\frac27}\tau^{\frac1{14}})(1-\zeta_7^3e^{\frac47}\tau^{\frac1{14}})}\Big)\\
&=\frac{1+\tau^{\frac3{14}}+3\tau^{\frac27}+3\tau^{\frac5{14}}+4\tau^{\frac37}+3\tau^{\frac12}+6\tau^{\frac47}+7\tau^{\frac9{14}}+6\tau^{\frac57}+3\tau^{\frac{11}{14}}+4\tau^{\frac67}+3\tau^{\frac{13}{14}}+3t+\tau^{\frac{15}{14}}+\tau^{\frac97}}{(1-\tau^{\frac12})^2}\end{align*}
and
\[\frac{1+3t^{\frac12}+3t}{(1-t^{\frac12})^2}=\VS{k\in\M}((k+1)+3k+3(k-1)^+)t^{\frac k2}=1+\VS{k\in\N}(7k-2)t^{\frac k2}.\]

The dimension formula says $\dim\cM(7)_k=14k-2$ for $k\in\N$.
\end{proof}

The dimension formula also says $\dim\cS(7)_k=\delta_{k,2}+14(k-\frac{12}7)-2)^+$ and
\[{\rm Dim}\,\cS(7)=\dfrac{3t^{\frac27}+3t^{\frac{11}{14}}+t^{\frac97}}{(1-t^{\frac12})^2}t^{\frac{12}7}\Big|_\M\]
thus
\[\natural:\cM(7)^*\eta^{\frac{24}7}|_\M\bto\cS(7).\]\\

At last, I conjecture
\[u_7^3u_7'=\VP{n\in\N}(1-q^n)^{-\frac37}\VP{n\in\M}(1-q^{7n+3})(1-q^{7n+4})(1-q^{7n+7}),\]
\[u_7'^3a_7=q^{\frac17}\VP{n\in\N}(1-q^n)^{-\frac37}\VP{n\in\M}(1-q^{7n+2})(1-q^{7n+5})(1-q^{7n+7}),\]
\[a_7^3u_7=q^{\frac37}\VP{n\in\N}(1-q^n)^{-\frac37}\VP{n\in\M}(1-q^{7n+1})(1-q^{7n+6})(1-q^{7n+7}).\]

\end{document}